\documentclass[12pt]{article}
\usepackage[utf8]{inputenc}
\usepackage{latexsym, amsmath, amsthm, amsfonts}
\usepackage{bbm}
\usepackage{enumerate, amssymb, xspace}
\usepackage[mathscr]{eucal}
\usepackage{graphicx}
\usepackage{rotating}
\usepackage{hyperref} %% Just avoid using dashes
\usepackage{color,colordvi}

\usepackage{mathtools}  %% for \xmapsto
  
\usepackage{dutchcal}   %% for lower case cal letters
  
\usepackage{tikz} 
\usepackage{braids}  %% NOTE: Must use \unloadcatcodes outside math mode before declaring a braid.
\usepackage{tikz-cd}
\usetikzlibrary{decorations.markings}
\usetikzlibrary{decorations.pathmorphing}
\usetikzlibrary{arrows.meta}  % for scaled arrows
  
\usepackage{pifont}
\usepackage{psfrag}
\usepackage{stackrel}
  
             \newcommand\Cite[2] {\cite[#1]{#2}}

\def\colorAlg   {black}
\def\colorMod   {red!80!black}

\def\colorMul   {red!70!black}

\def\colorRho   {yellow!80!black}
\def\symbMul    {circle}
\def\symbRho    {circle}

\def\widthMul   {0.12}
\def\widthRho   {0.14}
\def\widthObj   {1.5pt}  %  normal line thickness
\def\widthBcg   {4*\widthObj}  %  line thickness for "background"

\newcommand\drawMulLabelL[3]  {\fill[draw=\colorAlg,fill=\colorMul] (#1,#2)
                        \symbMul (\widthMul) node[left=-1pt,yshift=4pt] {$#3$}}
\newcommand\drawMulLabelR[3]  {\fill[draw=\colorAlg,fill=\colorMul] (#1,#2)
                        \symbMul (\widthMul) node[right,yshift=5pt] {$#3$}}

\newcommand\drawRhoLabel[3]  {\fill[draw=\colorAlg,fill=\colorRho] (#1,#2)
                        \symbRho (\widthRho) node[right=1pt,yshift=-1pt] {$#3$}}

\def\locscale   {1}

\newcommand\scopeArrow[2] {\begin{scope}[decoration={markings,mark=at position #1
                   with \arrow{#2}}]}  % put arrowhead of shape #2 at relative position #1

\newcommand{\leftturn}[1]{arc (0:90:#1)}
\newcommand{\leftcap}[1]{arc (0:180:#1)}

\newcommand{\rightcap}[1]{arc (180:0:#1)}

\newcommand{\leftmult}[1]{\leftcap{#1} arc (180:90:#1)}
\newcommand{\rightmult}[1]{\rightcap{#1} arc (0:90:#1)}

\def\erase{\draw[line width=3*\widthObj,white]}
\def\drawAlg{\draw[line width=\widthObj,\colorAlg]}
\def\drawMod{\draw[line width=\widthObj,\colorMod,behind path]}
\def\upto{to[out=90,in=270]}
\newcommand{\boxlabel}[1]{node[draw,fill=white]{$#1$} +(0.2,0.3)}

\def\AA            {{\ensuremath{\mathbb A}}}
\def\aA            {{\! A}}
\def\adj           {\alpha}  %% adjunction
\def\Ao            {A\opp_{\phantom i}}
\def\absi          {li\-ne\-ar-sim\-ple}

\def\ass           {\mathrm a}    % associator 
\def\be            {\begin{equation}}
\def\bea           {\be \begin{aligned}}
\def\bearl         {\begin{array}{l}}
\def\bearll        {\begin{array}{ll}}
\def\bearlll       {\begin{array}{lll}}
\def\bimod         {\text-\mathcal{bimod}}
\def\botIcoc       {{\circledast^{\coc}}}
\def\boticoc       {\,{\botIcoc}\,}
\def\br            {\beta}    % braiding - not c, as c is used for objects
\def\C             {{\ensuremath{\mathcal C}}}
\def\CA            {{\ensuremath{\mathcal C_{A}}}}   % all modules
\def\CAi           {{\ensuremath{\mathcal C_{\!A}}}}   % induced modules
\def\CatG          {\ensuremath{\ko\text-\mathcal C\hspace*{-1.6pt}at_{\!\Ga}}}

\def\cgmc          {$(\Ga\hspace*{-3pt},\hspace*{-1pt}\coc)$-gra\-ded-mo\-no\-i\-dal category} 
\def\cir           {\,{\circ}\,}

\def\coc           {\kappa}     % group 2-cocycle
\def\coco          {\omega}     % group 2-cochain
\def\Coco          {\Omega}     % group 2-cochain
\def\cocs          {\varpi}     % group 2-cocycle for stabilizer
\def\coe           {q_A^{}}   % coequalizer map
\def\Colon         {:\quad}

\def\comul         {\Delta}   % coproduct of a coalgebra
\def\comulA        {\comul^{\!\AA}} % coproduct of a comonad A
\def\comulb        {\overline\comul}   
\def\comulco       {\comul\ococ}
\def\comulM        {\boldsymbol\comul}
\def\comulS        {\comul^{\!\SS}} % coproduct of a comonad S
\def\comulT        {\comul^{\!\TT}} % coproduct of a comonad T
\def\Cong          {\,{\cong}\,}
\def\D             {{\ensuremath{\mathcal D}}}
\def\deltasep      {co\-pro\-duct-se\-pa\-ra\-ble}
\def\deltasepy     {co\-pro\-duct\text-se\-pa\-ra\-bi\-lity}
\def\df            {\,{:=}\,}

\def\dt            {\dot }
\def\dtm           {{\dt m}} 
\def\dtn           {{\dt n}}
\def\ee            {\end{equation}}
\def\eea           {\end{aligned} \ee}
\def\eear          {\end{array}}
\def\End           {{\mathrm{End}}}
\def\EndAA         {\End_\aA}
\def\Enumerate     {\renewcommand\leftmargini{1.34em} ~\\[-1.52em]
                   \begin{enumerate} \addtolength\itemsep{-6pt}}
\def\eps           {\varepsilon}
\def\Eps           {\eps^{\adj}}  % counit of adjunction
\def\epsA          {\eps^\AA} % counit of a comonad A
\def\epsS          {\eps^\SS} % counit of a comonad S
\def\epsT          {\eps^\TT} % counit of a comonad T
\def\eq            {\,{=}\,}
\newcommand\eqc[1] {\stackrel{\eqref{#1}}\cong}

\def\Eta           {\eta^{\adj}}  % unit of adjunction
\def\etaA          {\eta^\AA} % unit of a monad A
\def\etaT          {\eta^\TT} % unit of a monad T
\newcommand\FO[1]  {\mathrm F^{(#1\,#1\,#1)\,#1}_{\one,\one}}

\def\FxL           {\mathrm F_{\!L}}
\def\Ga            {\varGamma}
\newcommand\gr[1]  {{|#1|}}     % grade
\def\graco         {gra\-ded-com\-mu\-ta\-ti\-ve}
\def\Graco         {Gra\-ded-com\-mu\-ta\-ti\-ve}
\def\GSL           {Gra\-ded Schur Lemma}
\def\Hom           {{\mathrm{Hom}}}
\def\HomA          {\Hom_A}
\def\HomAA         {\Hom_\aA}
\def\HomC          {\Hom_\C}
\def\HomD          {{\ensuremath{\mathrm{Hom}_\D}}}
\def\I             {I\hspace*{-1.4pt}n\hspace*{-0.2pt}d}    % induction functor 
\def\id            {\text{\sl id}}
\def\Id            {\text{\sl Id}}
\def\ida           {\id_{\!A}}
\def\idaa          {\id_{\!\dt A}}

\newcommand\II[1]  {\I(#1)}
\newcommand\IIT[1] {\I_\TT(#1)}
\def\Im            {\mathrm{im}}
\def\iN            {\,{\in}\,}
\def\Ind           {\text-\Indu}
\def\inD           {\Indu\text-\hspace*{-0.9pt}}
\def\inDe          {\Indu^\circ\hspace*{-0.4pt}\text-\hspace*{-1.1pt}}
\def\Indu          {\mathcal i\hspace*{-1pt}\mathcal n\hspace*{-0.1pt}\mathcal d}
\def\inv           {^{-1}}
\def\Ker           {\mathrm{\ker}}
\def\ko            {{\ensuremath{\Bbbk}}}
\def\kotw          {\ko_\cocs}
\def\kox           {\ko^\times}
\def\lam           {\lambda}
\def\Mapsto        {\,{\xmapsto{~}}\,}
\def\Mod           {\text-\mathcal{mo}\hspace*{-0.2pt}\mathcal d}
\def\moD           {\mathcal{mod}\text-\hspace*{-0.9pt}}
\def\modGa         {\mathcal{mod}^{\circ}_{\!\Ga}\text-\hspace*{-0.9pt}}
\def\modeGa        {\mathcal{mod}_{\!\Ga}\text-\hspace*{-0.9pt}}
\def\mul           {\mu}    % product of an algebra
\def\mulA          {\mul^\AA} % product of a monad A
\def\mulb          {\overline\mul}
\def\mulM          {\boldsymbol\mul}
\def\mulco         {\mul\ococ}
\def\mulT          {\mul^\TT} % product of a monad T
\def\nE            {\,{\ne}\,}
\newcommand\nxl[1] {\\[#1mm]}

\def\obj           {\mathrm{Obj}}
\def\ocirc         {\overline\circ}
\def\ocir          {\,{\ocirc}\,}
\def\ococ          {^{[\coc]}}
\def\Ococ          {^{\![\coc]}}
\def\ococbr        {^{(\br;\coc)}}
\def\ococp         {^{[\coc']}}
\def\ohr           {\reflectbox{$\rho$}}

\def\ol            {\overline }
\def\one           {{\bf1}}
\def\onedim        {one-di\-men\-si\-o\-nal}

\def\opp           {^{\rm op}}
\def\ota           {\,\otA\,}
\def\otA           {{\otimes_{\!A}}}

\def\oti           {\,\otI\,}
\def\otI           {{\otimes}}
\def\otIcoc        {{\otimes^{\coc}}}
\def\oticoc        {\,{\otIcoc}\,}

\def\otk           {\,\otK\,}
\def\otK           {{\overline\otimes}} % tensor product on Kleisli cat  
\def\pl            {\,{+}\,}
\def\picc          {\mathrm{Pic}(\C)}
\def\Picc          {\mathcal Pic(\C)}
\def\pp            {\pi}    % generic idempotent
\def\qquand        {\qquad{\rm and}\qquad}
\newcommand\rarr[1]{\xrightarrow{~#1~}}
\newcommand\Rarr[1]{\,{\xrightarrow{\,#1\,}}\,}

\newcommand\RArr[1]{\,{\xRightarrow{\,#1\,}}\,}

\def\sjs           {$6j$-symbol}
\def\SS            {{\ensuremath{\mathbb S}}}
\def\sss           {\scriptscriptstyle}
\def\sta           {\varSigma}
\newcommand\stab[1]{\sta_{#1}}
\def\sZ            {s_{\sss(2)}}
\newcommand\SZ[2]  {\sZ^{#1,#2}}
\def\Times         {\,{\times}\,}

\def\toti          {\,\totI\,}
\def\totI          {{\widetilde\otimes}}
\def\TT            {{\ensuremath{\mathbb T}}}
\def\twi           {twisted}
\def\tZ            {t_{\sss(2)}}

\newcommand\TZ[2]  {\tZ^{#1,#2}}

\def\unlam         {\raisebox{0.65em}{\begin{turn}{180}{$\lambda$}\end{turn}}}
\def\unmul         {\raisebox{0.45em}{\begin{turn}{180}{$\mu$}\end{turn}}}
\def\unl           {\mathrm u_{\mathrm l}}         % left unitor
\newcommand\Unl[1] {\mathrm u_{{\mathrm l},#1}}
\def\unr           {\mathrm u_{\mathrm r}}         % right unitor
\newcommand\Unr[1] {\mathrm u_{{\mathrm r},#1}}
\def\Vect          {\text{vect}}
\def\VectG         {\ensuremath{\Vect_\Ga}}
\def\VectGe        {\ensuremath{\Vect_\Ga^{\scriptscriptstyle(0)}}}

\def\VectPPsi      {\ensuremath{\Vect_{\picc}^\Psi}}
\newcommand\void[1]{}

\def\wohr          {\ohr\Ococ}
\def\wrho          {\rho\ococ}

\def\zet           {{\mathbb Z}}
\def\zz            {{\ensuremath{\zet/2\zet}}}

\theoremstyle{definition}
\newtheorem{thm}{Theorem}[section]
\newtheorem{conv}[thm]{Convention}

\newtheorem{defi}[thm]{Definition}
\newtheorem{exa}[thm]{Example}
\newtheorem{lem}[thm]{Lemma}

\newtheorem{prop}[thm]{Proposition}
\newtheorem{rem}[thm]{Remark}

\definecolor{DarkViolet} {rgb}{0.580392,0.000000,0.827450}
\definecolor{ForestGreen}{rgb}{0.133333,0.545098,0.133333}
\definecolor{red3}       {rgb}{0.803921,0.000000,0.000000}

\newcounter{cmt}

\begin{document}

\numberwithin{equation}{section}
\numberwithin{thm}{section}

%%%%%%%%%%%%%%%%%%%%%%%%%%%%%%%%%%%%%%%%%%%%%%%%%%%%%%%%%%%%%%%%%%%%%%%%

~\vskip 3.3em

\begin{center}
\hspace*{-0.6em}{\bf \Large A Graded Schur Lemma and a graded-monoidal structure 
\nxl2
for induced modules over graded-commutative algebras}

\vskip 18mm

{\large \  \ J\"urgen Fuchs\,$^{\,a} \quad$ and $\quad$ Tobias Gr\o sfjeld\,$^{\,b}$ }

\vskip 12mm

 \it$^a$
 Teoretisk fysik, \ Karlstads Universitet\\
 Universitetsgatan 21, \ S\,--\,651\,88\, Karlstad
 \\[9pt]
 \it$^b$
 Matematiska institutionen, \ Stockholms Universitet\\
 S\,--\,106\,91\, Stockholm

\end{center}

\vskip 3.2em

\noindent{\sc Abstract}\\[3pt]
We consider algebras and Frobenius algebras, internal to a monoidal
category, that are graded over a finite abelian group. For the case that $A$ is a
twisted group algebra in a linear abelian monoidal category we obtain a graded 
generalization of the Schur Lemma for the category of induced $A$-modules.
We further show that if the monoidal category is braided and $A$ is commutative up 
to a bicharacter of the grading group, then the category of induced $A$-modules can
be endowed with a graded-monoidal structure that is twisted by the bicharacter.
In the particular case that the grading group is \zz, these findings reproduce
known results about superalgebras and super-monoidal structures.

%%%%%%%%%%%%%%%%%%%%%%%%%%%%%%%%%%%%%%%%%%%%%%%%%%%%%%%%%%%%%%%%%%%%%%%%
\newpage

\section{Introduction}

Algebras in monoidal categories and their representation theory have become
customary tools both in mathematics and in theoretical physics. While
for important aspects a Morita invariant approach via module categories over monoidal
categories is convenient, there are also issues which are naturally discussed
directly in terms of concrete representative algebras rather than their Morita
class. Specifically, a significant property an algebra may possess that is not 
Morita invariant is commutativity.
  
The category of right modules over a commutative algebra $A$ in a (suitably nice)
\ko-linear braided monoidal category \C\ inherits a natural monoidal structure from 
the tensor product over $A$ of $A$-bimodules. In this note we study how this feature
generalizes to the case of algebras $A$ that are \emph{\graco} in the sense that
$A \eq \bigoplus_{i\in\Ga\!} A_i$ is graded over a finite abelian group $\Ga$ such
that the multiplication on $A$ is grade-preserving and (braided-)commutative up to
a normalized two-cocycle on $\Ga$ 
(see Definitions \ref{def:Ga-graded.alg} and \ref{def:graco}).
We first consider the non-full subcategory $\modGa A$ of the category $\modeGa A$
of $\Ga$-graded right $A$-modules that is obtained by admitting as morphisms only 
grade-preserving
module morphisms. We show in Proposition \ref{prop:tildeoti} that if the tensor
product $\otA$ of $A$-bimodules exists and $A$ is \graco\ with respect to a cocycle
$\coc$ that is a bicharacter, then $\otA$ furnishes a monoidal structure on
$\modGa A$.

\medskip

Our basic goal in this note is then to obtain similar results for categories of
$A$-modules on which no restriction to grade-preserving morphisms is imposed. To
achieve this we consider instead of the category $\modeGa A$ of all right $A$-modules
its full subcategory $\inD A$ of \emph{induced} $A$-mo\-dules.
Our first main result (Proposition \ref{prop:dimEndAA=ko.stab} and Theorem 
\ref{thm:GSL}) then concerns the structure of the morphism spaces in such categories
of induced modules, without reference to a braiding on \C: For any twisted group
algebra $A$ (as defined in  Example \ref{exa:twistedgroupalg})
in a \ko-li\-ne\-ar abelian monoidal category \C\ with simple monoidal unit
the space of module morphisms between any two simple induced $A$-modules $\II x$
and $\II y$ is either zero or isomorphic as a $\Ga$-graded \ko-vector space to
the endomorphims of $\II x$ (as well as to those of $\II y$) which, in turn, is 
isomorphic as a $\stab x$-graded \ko-algebra to the group algebra $\ko\, \stab x$ 
twisted by a two-cocycle, where $\stab x$ is the subgroup
$\{ i \iN \Ga \,|\, x \oti A_i \Cong x \}$ of $\Ga$. Moreover, any non-zero 
homogeneous module morphism between simple induced $A$-modules is an isomorphism.
This result constitutes a generalization of the Schur Lemma for abelian categories,
and accordingly we refer to it as the \emph{\GSL}. As we point out in Remark 
\ref{rem:SSL}, in the case $\Ga\eq\zz$ this comprises several variants of what is 
known as the Super Schur Lemma, which is e.g.\ relevant in the study of fermionic
topological phases in condensed matter physics.

\medskip

For a generic \graco\ algebra $A$, it is natural to ask whether the tensor product
$\totI$ that the category $\inDe A$ of induced $A$-modules with grade-preserving
module morphisms inherits as a full subcategory of $\modGa A$
can be extended to all of $\inD A$. This turns out to be possible only if $A$
is actually commutative (see Lemma \ref{lem:monmonad4commFrob}).
However, this observation can be upgraded to the \graco\ situation, thus leading to
our second main result (Theorem \ref{thm:II}) which states: For an algebra $A$ in a 
braided monoidal category that is \graco\ with respect to a bicharacter $\coc$ of
a finite abelian group $\Ga$, the category $\inD A$ can be endowed with the 
structure of a \emph{\cgmc}. Furthermore, 
if the tensor product of $A$-bimodules exists and $A$ is \deltasep\ Frobenius, 
then restriction to the subcategory $\inDe A$ reproduces the tensor product $\totI$.

Here the relevant notion of a linear \cgmc\ is the following. Such a category is,
by definition, a category enriched over the category $\VectGe$ of $\Ga$-graded vector
spaces with grade-preserving linear maps, in which the morphism spaces are graded and
the tensor product of morphisms is twisted by the bicharacter $\coc$ \cite{puty,brDa}.
For the details see Definitions and \ref{def:cocgracat} and \ref{def:cocgramoncat}.
The notion of \cgmc\ is thus in particular substantially different from other
structures of graded category that appear in the literature, both from the one
in which the objects are graded (see e.g.\ \Cite{Ch.\,9.13}{EGno} for the case of 
fusion categories) and which e.g.\ underlies the notion of a $\Ga$-equivariant 
category, and from the ``$\Ga$-graded monoidal structure'' as defined in \cite{frwa} 
and studied in e.g.\ \cite{cegO}. Instead, it is a generalization from $\zz$ to an 
arbitrary finite abelian group of the notion of a supermonoidal structure, which arises
in the context of monoidal supercategories of the form considered in \cite{brEl}.

The additional assumption of $A$ being Frobenius as well as \deltasep,
in the sense of Definitions \eqref{def:frob} and \eqref{def:sepa}, 
that is made in the second part of Theorem \ref{thm:II} is in particular
satisfied for twisted group algebras. Technically, its relevance derives from the 
fact that the tensor product of $A$-bimodules, if it exists, can be realized as the
image of an idempotent, see Lemma \ref{lem:mAn=Im(p)}. 
Note that \deltasep\ Frobenius algebras in monoidal categories play a crucial role in
various different contexts, such as in the study of weak monoidal Morita equivalence 
of tensor categories \cite{muge8} and of autonomous pseudomonoids in monoidal 
bicategories \cite{stre8}, and in two-dimensional conformal quantum field theory 
\cite{fuRs4,henri}.

\bigskip 

The rest of this paper is organized as follows.
In Sections \ref{sec:prel-1}\,--\,\ref{sec:prel-3} we present pertinent background
information on graded algebras and graded Frobenius algebras in monoidal categories,
on their categories of graded modules and on the tensor product of their bimodules.
This includes a few illustrating examples, in particular the twisted group algebras
in Examples \ref{exa:twistedgroupalg} and \ref{exa:twistedgroupalg2}, as well as
a discussion of analogues for modules (Lemma \ref{lem:mtilde})
and for Frobenius algebras (Lemma \ref{def:lem:twist-Frobalg}) of the 
cocycle-twisting of algebras.
We also collect, in Convention \ref{conv}, some conventions and notations
that are used throughout the text.
Section \ref{sec:prel-4} presents the notion of a \graco\ algebra and a few results
about them; e.g.\ a twisted group algebra is \graco\ iff the twisting cochain is 
a bicharacter (Lemma \ref{lem:graco-iff-bichar}).
Sections \ref{sec:prel-5} and \ref{gr-enriched} provide supplementary background on
monoidal structures on monads and on the appropriate notion of graded-monoidal 
category, respectively.

In Section \ref{sec:oti} we exhibit the tensor product on the category $\modGa A$ 
for a \graco\ algebra $A$, including the description as the image of an idempotent
(Proposition \ref{prop:tildeoti}(ii)) in case $A$ is \deltasep\ Frobenius.
Section \ref{sec:inDA} discusses aspects of the category $\inD A$ of induced
$A$-modules for the case that $A$ is a twisted group algebra, culminating in the
\GSL, Theorem \ref{thm:GSL}.

Section \ref{sec:gil} is devoted to establishing Theorem \ref{thm:II}, i.e.\
providing the structure of a \cgmc, in the sense of 
Definition \ref{def:cocgramoncat}, on the category of induced modules over an algebra
$A$ that is \graco\ with respect to a bicharacter $\coc$. This is achieved by realizing
the category $\inD A$ as the Kleisli category over the monad $\AA$ with underlying
endofunctor ${-}\oti A$, and identifying, in Proposition \ref{prop:twistedinterchange},
a corresponding graded variant of a monoidal
structure on $\AA$ that implies a $\coc$-twisted interchange law.

     \newpage  % XXX 
%%%%%%%%%%%%%%%%%%%%%%%%%%%%%%%%%%%%%%%%%%%%%%%%%%%%%%%%%%%%%%%%%%%%%%%%

\section{Preliminaries} \label{sec:prel}

\subsection{Graded algebras} \label{sec:prel-1}

Let $\C \eq (\C,\otimes,\one,\ass,\unl,\unr)$ be an additive monoidal category.
Recall that a (unital, associative) algebra in \C\ is a triple 
$A \eq (\dt A,\mul,\eta)$ consisting of an object $\dt A\iN\C$, a multiplication
morphism $\mul \iN \Hom(\dt A\oti \dt A,\dt A)$ and a unit morphism
$\eta \iN \Hom(\one,\dt A)$ satisfying 
$\mul \cir (\mul \oti \idaa) \eq \mul \cir (\idaa \oti \mul)
    $\linebreak[0]$ {\circ}\,
    \ass_{\dt A,\dt A,\dt A}$
(associativity) and $\mul \cir (\eta \oti \idaa) \cir \Unl{\dt A}\inv \eq \idaa 
\eq \mul \cir (\idaa \oti \eta) \cir \Unr{\dt A}\inv$ (unit properties).

Below we will adhere to the following conventions, unless stated otherwise:

\begin{conv} \label{conv}
 \Enumerate
 \item
Throughout this paper, \ko\ is an algebraically closed field.
A generic category will be assumed to be additive and \ko-linear.
 \item
Whenever no confusion can arise, we suppress the associator $\ass$ 
and unitors $\unl$ and $\unr$ of a monoidal category. 
 \item
When $x$ is an object with additional structure, such as an algebra or module, we 
denote the underlying plain object by $\dt x$. However, we often abuse notation by
denoting the object underlying an algebra $A \eq (\dt A,\mul,\eta)$ just by the symbol
$A$, and similarly for modules.
 \item \label{conv-i4}
For any \ko-linear monoidal category we assume that $\Hom(\one,\one) \,{\cong}\, \ko$,
and for any algebra we assume that its underlying object is non-zero.
 \item
Throughout this paper, $\Ga$ is a 
finite abelian group and $\coc$ is a normalized two-cocycle on $\Ga$. We write 
the group product of $\Ga$ as addition; in particular the unit element is 
denoted by $0 \iN \Ga$. If $\Ga$ occurs in a \ko-linear setting, we assume 
that the characteristic of \ko\ and the order of $\Ga$ are coprime. 
 \item
In Section \ref{sec:inDA} $A$ is a twisted group algebra in the sense of Example
\ref{exa:twistedgroupalg}, while in Sections \ref{sec:oti} and \ref{sec:gil} $A$ is
\graco\ in the sense of Definition \ref{def:graco} and it is assumed that
the tensor product $\otA$ of $A$-bimodules exists.
 \end{enumerate}
\end{conv}

For an object $c\iN\C$ with a direct sum decomposition $c \eq \bigoplus_{i\!} c_i$,
we denote by $e_i^c$ and $r^i_c$ the embedding and restriction morphisms for
$c_i$ as a retract of $c$ and set $e_i^c \cir r^i_c \,{=:}\, p_i^c$. By definition 
we have $r^i_c \cir e_j^c \eq \delta_{i,j}\,\id_{c_i}$, so that the morphisms $p_i^c$
are mutually orthogonal idempotents, $p_i^c \cir p_j^c \eq \delta_{i,j}\, p_i^c$\,.
If the direct summands are labeled by the group $\Ga$, i.e.\
  \be
  c = \bigoplus_{i\in\Ga\!} c_i \,,
  \label{eq:def:gradedobject}
  \ee
we call $c$ a \emph{$\Ga$-graded object}.
The tensor product of two $\Ga$-graded objects $c$ and $d$ is naturally $\Ga$-graded,
with $(c\oti d)_i \eq \bigoplus_{j\in\Ga\!} c_j \oti d_{i-j}$.
A morphism $f\colon c\Rarr~d$ between $\Ga$-graded objects is \emph{grade-preserving}
iff it satisfies
  \be
  r^j_d \circ f \circ e^c_i = 0  \quad \text{for all}~ i,j\iN\Ga ~\text{with}~
  i \,{\ne}\, j \,.
  \label{eq:def:evenmor}
  \ee
Borrowing terminology from the case $\Ga \eq \zz$, we call a grade-preserving
morphism synonymously also a \emph{homogeneously even} morphism.

\begin{defi} \label{def:Ga-graded.alg}
A \emph{$\Ga$-graded algebra} in a monoidal category \C\ is an algebra
$A \eq (\dt A,\mul,\eta)$ in \C\ for which the underlying object $\dt A\iN\C$ 
is $\Ga$-graded and the multiplication morphism $\mul$ is homogeneously even.
\end{defi}

Thus a $\Ga$-graded algebra has a direct sum decomposition 
$A \eq \bigoplus_{i\in\Ga\!} A_i$. When no confusion can arise, we abbreviate
$e_i^A \,{=:}\, e_i$, $r^i_\aA \,{=:}\, r^i$ and $p_i^A \,{=:}\, p_i$.
That $\mul$ is homogeneously even then means that the morphisms
  \be
  \mul_{i,j}^k := r^k \cir \mul \cir (e_i \oti e_j) \,\in \Hom(A_i\oti A_j,A_k)
  \label{eq:mul-ijk}
  \ee
satisfy $\mul_{i,j}^k \eq 0$ if $k \nE i \pl j$, for $i,j,k\iN\Ga$.
Accordingly we will use the abbreviation
  \be
  \mul_{i,j}^{i+j} \,{=:}\, \mul_{i,j} \,.
  \label{eq:abbr.mulij}
  \ee
It follows from the unit axioms that $e_0 \cir r^0 \cir \eta \eq \eta$
and thus, with the natural convention that the monoidal unit is $\Ga$-graded with
$\one_j \eq \delta_{j,0}\, \one$, that $\eta$ is automatically homogeneously even.

\begin{rem} \label{rem:twist-alg}
Let $A \eq (\dt A,\mul,\eta)$ be a $\Ga$-graded algebra in a \ko-linear
monoidal category \C, and let $\coc\colon \Ga\Times\Ga \Rarr~ \kox$ be a normalized
two-cocycle on $\Ga$. Then as a direct generalization of a standard result for 
\ko-algebras (compare e.g.\ Chapter 2 of \cite{KArp2} for the case of group algebras),
the triple $A\ococ \df (\dt A,\mulco,\eta)$ with cocycle-twisted multiplication
  \be
  \mulco_{i,j} := \coc(i,j)\, \mul_{i,j}
  \label{eq:def:wmulijk}
  \ee
defines a structure of an associative algebra in \C. To see that $\mulco$ is an 
associative multiplication on $\dt A$, use
$\mul \eq \sum_{i,j\in\Ga} e_{i+j} \cir \mul_{i,j} \cir (r^i \oti r^j)$ and the
analogous decomposition of $\mulco$ to write
  \be
  \begin{aligned}
  \mulco \circ (\mulco \oti \ida)
  & = \!\sum_{i,j,k\in\Ga}\! \coc(i,j)\, \coc(i{+}j,k)\, e_{i+j+k} \circ
  [ \mul_{i+j,k} \cir (\mul_{i,j} \oti \id_{A_k}) ] \circ (r^i \oti r^j \oti r^k) \,,
  \\
  \mulco \circ (\ida \oti \mulco)
  & = \!\sum_{i,j,k\in\Ga}\! \coc(j,k)\, \coc(i,j{+}k)\, e_{i+j+k} \circ
  [ \mul_{i,j+k} \cir (\id_{A_i} \oti \mul_{j,k}) ] \circ (r^i \oti r^j \oti r^k)
  \,.
  \end{aligned}
  \ee
By associativity of $\mul$ we have $\mul_{i+j,k} \cir (\mul_{i,j} \oti \id_{A_k}) 
\eq \mul_{i,j+k} \cir (\id_{A_i} \oti \mul_{j,k})$, while by the closedness
$\mathrm d\coc \eq 1$ of $\coc$ the cocycle factors in the two expressions are equal.
That $\eta$ is a unit for $\mulco$ follows from
$\eta \eq p_0 \cir \eta$ together with the fact that $\coc$ is normalized,
i.e.\ $\coc(i,0) \eq 1 \eq \coc(0,i)$ for $i\iN\Ga$.
\end{rem}

\begin{lem} \label{lem:cohomologous}
Let $A \eq (\dt A,\mul,\eta)$ be a $\Ga$-graded algebra in a \ko-linear monoidal
category \C\ and let $\coc$ and $\coc'$ be normalized cocycles on $\Ga$.
If $\coc$ and $\coc'$ are cohomologous, then the algebras $A\ococ$ and $A\ococp$
obtained as in Remark \ref{rem:twist-alg} are isomorphic
by a grade-preserving morphism.
Conversely, provided that the morphism $\mul_{i,j}$ is non-zero for all $i,j\iN\Ga$,
if $A\ococ$ and $A\ococp$ are isomorphic as algebras by a grade-preserving morphism,
then $\coc$ and $\coc'$ are cohomologous.
\end{lem}
    
\begin{proof}
That $\coc$ and $\coc'$ are cohomologous means that
$\coc'(i,j) \eq \tau(i)\, \tau(j)\, \tau\inv(i{+}j)\, \coc(i,j)$
for some function $\tau\colon \Ga\Rarr~\kox$ with $\tau(0)\eq1$. Thus the 
multiplications on $A\ococ$ and on $A\ococp$ are related by
$\mul\ococp_{i,j} \eq \tau(i)\, \tau(j)\, \tau\inv(i{+}j)\, \mul\ococ_{i,j}$. It 
follows that the
automorphism $\varphi \,{:=}\, \bigoplus_{i\in\Ga} \tau\inv(i)\, \id_{A_i}$ of the
object $\dt A$ satisfies $\big( \mul\ococp \cir (\varphi\oti\varphi) \big)_{i,j}
\eq \big( \varphi \cir \mul\ococ \big)_{i,j}$
for all $i,j\iN\Ga$ and is thus a morphism of associative algebras from $A\ococ$ to 
$A\ococp$. Because of $\tau(0)\eq1$, $\varphi$ is also compatible with the unit.
To see the converse, assume that $A\ococ$ and $A\ococp$ are isomorphic as objects
by a grade-preserving morphism $\varphi$. Owing to the equality
$\mul\ococp_{i,j} \eq \coc'(i,j)\, \coc(i,j)\inv \mul\ococ_{i,j}$, in order that 
$\varphi$ is an algebra morphism it must act as a scalar on the image of $\mu_{i,j}$.
If $\mul_{i,j} \,{\ne}\, 0$ for all $i,j\iN\Ga$, we obtain in this way a one-cochain
by which $\coc'$ and $\coc$ are cohomologous.
\end{proof}

\begin{exa} \label{exa:subgroupgrading}
For any normal subgroup $H$ of a finite group $G$, the group algebra $\ko[G]$ is 
a $G/H$-graded algebra in the category of \ko-vector spaces:
$\ko[G] \eq \bigoplus_{gH\in G/H} \ko[gH]$.
\end{exa}

\begin{exa} \label{exa:twistedgroupalg}
A $\Ga$-graded algebra in the category \Vect\ of finite-dimensional \ko-vector 
spaces can be seen as an algebra in the category \VectGe\ of finite-dimensional 
$\Ga$-graded \ko-vector spaces and grade-preserving \ko-linear maps. The following
generalization (compare \Cite{Example\,7.8.3}{EGno}) is our main, and motivating,
class of examples. 
Denote by $\picc$ the Picard group of \C, by which we mean the group of isomorphism
classes of invertible objects of \C\ (not to be confused with the group of equivalence
classes of invertible \C-module categories, which is also called Picard group and
is denoted by {\sf Pic}(\C) in e.g.\ \Cite{8.27.2}{EGno}).
Consider an object $A$ that is pointed, i.e.\ a direct sum
  \be
  A_{(\Ga)} := \bigoplus_{i\in\Ga\!} L_i
  \ee
of invertible objects $L_i$, with $\Ga$ a finite subgroup of $\picc$,
and require that $L_{i+j} \Cong L_i \oti L_j$ and $L_i \,{\not\cong}\, L_j$ for
$i \nE j$. (Note that, since according to Convention \ref{conv}.\ref{conv-i4} the
monoidal unit $\one$ is \absi, i.e.\ with endomorphism space being isomorphic to
\ko, every invertible object is \absi\ as well.)
 \\[2pt]
The Picard subcategory $\Picc$, i.e.\ the full monoidal subcategory of \C\ consisting
of direct sums of invertible objects, is equivalent as an abelian category to the 
category \VectPPsi\ of $\picc$-graded \ko-vector spaces with an associator that 
defines a three-cocycle $\Psi \iN Z^3(\picc,\kox)$; monoidally equivalent categories
give rise to cohomologous three-cocycles. Assume further that there is a two-cochain 
$\coco \iN C^2(\Ga,\kox)$ 
on $\Ga$ that satisfies $\mathrm d\coco \eq \Psi|_\Ga^{}$, i.e.\ in components, 
  \be
  \coco(i,j)\,\coco(i\pl j,k) = \coco(j,k)\, \coco(i,j\pl k)\, 
  \Psi(i,j,k) \qquad \text{for}~~ i,j,k\iN\Ga . 
  \label{eq:def:2coc}
  \ee
Now select for each pair of elements $i,j\iN\Ga$ a basis $\lam_{i,j}$
of the one-dimensional morphism space $\Hom(L_i\oti L_j,L_{i+j})$. Then setting 
  \be
  \mul_{i,j} := \coco(i,j)\, \lam_{i,j}
  \label{eq:mul=coco.lam}
  \ee
for $i,j\iN\Ga$ defines a structure of a $\Ga$-graded algebra on the 
$\Ga$-graded object $A_{(\Ga)} \iN \Picc$ (associativity precisely amounts to 
\eqref{eq:def:2coc}). When defining $\eta$ as the embedding of the subobject 
$L_0 \Cong \one$ into $A_{(\Ga)}$, this multiplication is also unital.
We refer to such a $\Ga$-graded algebra as a \emph{twisted group algebra}.
 \\[2pt]
For $\coc$ any normalized two-cocycle on $\Ga$, owing to 
$\mathrm d(\coco\,\coc) \eq \mathrm d\coco\, \mathrm d\coc \eq \mathrm d\coco$,
$A_{(\Ga)}$ endowed with the $\coc$-twisted multiplication $\mul\ococ$ \eqref{eq:def:wmulijk}
is again a twisted group algebra. In fact, every semisimple 
indecomposable algebra in \VectPPsi\ is such an algebra \Cite{Ex.\,2.1}{ostr5}.
Replacing $\mul$ by $\mul\ococ$ amounts to changing the chosen bases $\lam_{i,j}$ to
$\coc(i,j)\,\lam_{i,j}$.  Also, if $\coc \eq \mathrm d\chi$ is a coboundary, then 
by Lemma \ref{lem:cohomologous} the twisting by $\coc$ can be absorbed into an 
automorphism $\sum_{i\in\Ga}\chi(i)\,\id_{\!L_i}$ of $A_{(\Ga)}$, so that 
$(A_{(\Ga)},\mul)$ and $(A_{(\Ga)},\mul^{[\mathrm d\chi]})$ are isomorphic as algebras.
 \\[2pt]
Note that this works in fact also for non-abelian groups $A_{(\Ga)}$. Later on (see
Example \ref{exa:twistedgroupalg3}) we will consider braided categories \C, for
which $\picc$, and thus $\Ga \,{\le}\, \picc$, is abelian.
\end{exa}

\begin{exa} \label{exa:zz}
Consider the special case $\Ga \eq \picc \eq \zz$ of Example \ref{exa:twistedgroupalg},
and let $L \iN \C$ be the (up to isomorphism unique) order-2 invertible object.
Then the class of $\Psi$ in $H^3(\zz,\kox)$ is determined by the value
$\FO L \iN \{\pm1\}$ of the gauge-independent \sjs\ $\FO L$ of the
semisimple category $\Picc$.\,%
 \footnote{~We use standard conventions for the \sjs s 
 of a semisimple monoidal category, see e.g.\ Appendix A of \cite{fuGr}.}
The object $A \eq \one \oplus L$ carries an algebra structure iff $\FO L \eq 1$, 
and in this case the algebra is of the form 
described in Example \ref{exa:twistedgroupalg} with $\Psi$ and $\omega$ trivial.
As a concrete realization, take the semisimple modular tensor
category $\C(\mathrm A_1^{(1)},k)$ whose objects correspond to the integrable
representations of the untwisted affine Lie algebra $\mathrm A_1^{(1)}$ at 
non-negative integral level $k$ (compare e.g.\ \Cite{Ex.\,8.18.5\,\&\,8.18.6}{EGno}). 
This category has $k{+}1$ non-isomorphic simple objects $U_\nu$, labeled by 
$\nu \iN \{0,1,...\,,k\}$, has Picard group $\zz$ with $L \eq U_k$, and
$\FO L \eq (-1)^k_{}$. Thus the object $A \eq \one \oplus L$ is an algebra in
$\C(\mathrm A_1^{(1)},k)$ iff $k$ is even.
\end{exa}

\begin{exa}
A category \C\ may admit a \zz-graded algebra with underlying object of the form 
$A \eq A_0 \oplus A_1$ such that $A_1$ is simple while 
$A_0 \,{\cong}\, \one \oplus L$ with $L$ an order-2 invertible object and 
$L \oti A_1 \,{\cong}\, A_1 \,{\cong}\, A_1 \oti L$. A category for which this 
happens is, in the notation of Example \ref{exa:zz}, $\C(\mathrm A_1^{(1)},4)$, with 
$L \eq U_4$ and $A_1 \eq U_2$ \Cite{Sect.\,3.6.2}{fuRs4}. In this case one
has $U_2 \oti U_2 \,{\cong}\, U_0 \,{\oplus}\, U_2 \,{\oplus}\, U_4$
and thus $A \,{\cong}\, A_1 \oti A_1$, so that
in particular $A$ is Morita equivalent to the
trivial algebra $\one \iN \C(\mathrm A_1^{(1)},4)$.  We do not know of a category
with a \zz-graded algebra of this type that is not Morita equivalent to $\one$.
\end{exa}

\begin{conv} \label{conv:mulmatrix}
If for a $\Ga$-graded algebra $A \eq \bigoplus_{i\!} A_i$ the morphism spaces 
$\Hom(A_i\oti A_j,A_{i+j})$ are zero- or \onedim\ for all $i,j\iN\Ga$, then,
given any choice $\{ \lam_{i,j} \}$ of bases of these spaces, we set
$\mul_{i,j}^{} \,{\equiv}\, \mul_{i,j}^{i+j} \,{=:}\, \mulb_{i,j} \, \lam_{i,j}$
with numbers $\mulb_{i,j} \iN \ko$ and call the matrix
$\mulM \,{:=}\, \big( \mulb_{i,j} \big)_{\!i,j\in\Ga}$
the \emph{matrix representation} of $\mu$. By abuse of notation we use
the same symbols and terminology whenever $\mul_{i,j} \eq 0$, independently 
of the dimensionality of $\Hom(A_i\oti A_j,A_{i+j})$.
\end{conv}

\begin{exa} \label{exa:degenerate1}
Let $x\iN\C$ be arbitrary. Then the object $X \,{:=}\, \one \oplus x$
endowed with the multiplication $\mul$ with matrix representation
  \be
  \mulM = \begin{pmatrix} 1 & 1 \\ 1 & 0 \end{pmatrix}
  \label{eq:mulM.degenerate}
  \ee
is a unital associative algebra, and it is \zz-graded with $X_0 \eq \one$ and 
$X_1 \eq x$. This algebra $X$ can be seen as an analogue of the quotient ring 
$\ko[t]/{<}t^2{>}$; likewise there are analogues of 
$\ko[t]/{<}t^\ell{>}$ for any integer $\ell \,{>}\, 2$.
\end{exa}

Recall that a right module $m$ over an algebra $A$ in \C\ consists of an object 
$\dtm\iN\C$ and a representation morphism $\ohr\colon \dtm \oti A\Rarr~ \dtm$ 
satisfying $\ohr \cir (\ohr \oti \ida) \eq \ohr \cir (\id_\dtm \oti \mul)$ and
$\ohr \cir (\id_\dtm \oti \eta) \eq \id_\dtm$. Analogously, for a left $A$-module
$n \eq (\dtn,\rho)$ the representation morphism $\rho\colon A\oti \dtn\Rarr~ \dtn$
obeys $\rho \cir (\ida \oti \rho) \eq \rho \cir (\mul \oti \id_\dtn)$ and
$\rho \cir (\eta \oti \id_\dtn) \eq \id_\dtn$. A bimodule is both a left and a
right module, with commuting left and right actions.
A morphism $f\colon m \Rarr~ m'$ of right $A$-modules is a morphism of the 
underlying objects that commutes with the $A$-action, i.e.\ $f\colon \dtm\Rarr~ \dtm'$
such that $\ohr_n \cir (f \oti \ida) \eq f \cir \ohr_m$; we denote the set of such
module morphisms by $\HomAA(m,n)$.
We will often refer to a right $A$-module just as an $A$-module.

\begin{defi} \label{def:Ga-module}
Let $A$ be a $\Ga$-graded algebra in a monoidal category \C. A \emph{$\Ga$-graded
right module} over $A$ is a right $A$-module $(\dtm,\ohr)$ for which the object 
$\dtm\iN\C$ is $\Ga$-graded (i.e.\ has a direct sum decomposition as in 
\eqref{eq:def:gradedobject}) and the action morphism $\ohr$ is homogeneously even.
The \emph{category $\modGa A$ of $\Ga$-graded right $A$-modules} is the category
that has $\Ga$-graded right $A$-modules as objects and homogeneously even
$A$-module morphisms as morphisms.
\end{defi}

More explicitly, a $\Ga$-graded right module $(\dtm,\ohr)$ has a direct sum 
decomposition $\dtm \eq \bigoplus_{i\in\Ga\!} m_i$ such that the morphisms
$\ohr_{i,j}^k \df r^k_m \cir \ohr \cir (e_i^m \oti e_j^{}) \iN \Hom(m_i \oti A_j,m_k)$
satisfy $\ohr_{i,j}^k \eq 0$ if $k \nE i \pl j$, for $i,j,k\iN\Ga$. And a
homogeneously even $A$-module morphism is a module morphism $f$ between
$\Ga$-graded right $A$-modules $(\dtm,\ohr)$ and $(\dtm',\ohr')$ such that
$r^j_{m'} \cir f \cir e^m_i \eq 0$ if $j\nE i$, for $i,j\iN\Ga$.
Similarly as done for graded algebras in \eqref{eq:abbr.mulij}, for a $\Ga$-graded 
right module $(c,\ohr)$ we abbreviate $\ohr_{i,j}^{i+j} \,{=:}\, \ohr_{i,j}$.

The notion of cocycle-twisting of algebras in Remark \ref{rem:twist-alg} has a
counterpart for modules. The two constructions are related as follows:

\begin{lem} \label{lem:mtilde}
Let $m \eq (\dtm,\ohr)$ be a $\Ga$-graded right module over a $\Ga$-graded algebra
$A$ in a \ko-linear monoidal category \C, and let $\coc$ be a normalized 
two-cocycle on $\Ga$. The $\coc$-twisted action
  \be
  \wohr := \sum_{i,j\in\Ga} e_{i+j}^{\dtm} \circ \wohr_{i,j}
  \circ (r^i_\dtm \oti r^j_{}) \qquad \text{with} \qquad
  \wohr_{i,j} := \coc(i,j)\, \ohr_{i,j}
  \label{eq:def:wrho}
  \ee
endows the object $\dtm\iN\C$ with the structure of a right $A\ococ$-module
  \be
  m\ococ = (\dtm,\wohr) \,.
  \label{eq:mtilde}
  \ee
Analogously, for $n \eq (\dtn,\rho)$ a $\Ga$-graded left $A$-module, the object 
$\dtn$ comes with the structure of a left $A\ococ$-module $n\ococ \eq (\dtn,\wrho)$
with $\coc$-twisted left action
$\wrho \df \sum_{i,j\in\Ga} e_{i+j}^\dtn \cir \wrho_{i,j} \cir (r^i_{} \oti r^j_\dtn)$ 
with $\wrho_{i,j} \df \coc(i,j)\, \rho_{i,j}$.
\end{lem}

\begin{proof}
The sequence
  \be
  \begin{aligned}
  \wohr \cir (\wohr \oti \ida)
  & = \sum_{i,j,k\in\Ga}\! \coc(i{+}j,k)\,\coc(i,j) \, e^\dtm_{i+j+k} \circ
  \ohr_{i+j,k} \circ (\ohr_{i,j}\oti  \id_{A_k} ) \circ (r^i_\dtm \oti r^j \oti r^k)
  \\
  & = \sum_{i,j,k\in\Ga}\! \coc(i,j{+}k)\,\coc(j,k) \, e^\dtm_{i+j+k} \circ
  \ohr_{i,j+k} \circ ( \id_{\dtm_i} \oti \mul_{j,k}) \circ (r^i_\dtm \oti r^j \oti r^k)
  \\
  & = \sum_{i,j,k\in\Ga}\! e^\dtm_{i+j+k} \circ
  \wohr_{i,j+k} \circ (\id_{m_i} \oti \mulco_{j,k}) \circ (r^i_\dtm \oti r^j \oti r^k)
  \\
  & = \wohr  \cir (\id_\dtm \oti \mulco)
  \end{aligned}
  \ee
of equalities shows that $\wohr$ is compatible with the twisted multiplication 
$\mulco$. Here the third equality uses the cocycle condition for $\coc$ as well as the
compatibility of $\ohr$ with $\mul$. The compatibility of $\wohr$ with the unit 
is seen similarly by using that the cocycle $\coc$ is normalized.
The calculation for the left module case is analogous.
\end{proof}

%%%%%%%%%%%%%%%%%%%%%%%%%%%%%%%%%%%%%%%%%%%%%%%%%%%%%%%%%%%%%%%%%%%%%%%%

\subsection{Frobenius algebras and separability} \label{sec:prel-2}

A (counital, coassociative) coalgebra in \C\ is a triple 
$B \eq (B,\comul,\eps)$ consisting of an object $B \iN \C$ and morphisms 
$\comul \iN \Hom(B,B\oti B)$ and $\eps \iN \Hom(B,\one)$ satisfying
$(\comul \oti \ida) \cir \comul \eq (\ida \oti \comul) \cir \comul$ and
$(\eps \oti \ida) \cir \comul \eq \ida \eq (\ida \oti \eps) \cir \comul$.

\begin{defi} \label{def:frob}
A \emph{Frobenius algebra} in \C\ is a quintuple $F \eq (F,\mul,\eta,\comul,\eps)$
such that $(F,\mul,\eta)$ is an algebra, $(F,\comul,\eps)$ is a coalgebra, and
  \be
  (\mul \oti \ida)  \circ (\ida \oti \comul) = \comul \circ \mul
  = (\ida \oti \mul) \circ (\comul \oti \ida) \,,
  \label{defi:frob}
  \ee
i.e.\ the coproduct is a morphism of $F$-bimodules.
\end{defi}

\begin{defi} \label{def:sepa}
A \emph{separable} algebra in \C\ is an algebra $(A,\mul,\eta)$ such that the
multiplication $\mul$ splits as a morphism of $A$-bimodules, i.e.\ there is a 
bimodule morphism $\unmul \iN \Hom(A,A\oti A)$ such that $\mul \cir \unmul \eq \ida$.
A \emph{\deltasep} algebra, also called \emph{$\Delta$-separable} algebra, in \C\ is 
an algebra-coalgebra $(A,\mul,\eta,\comul,\eps)$ which as an algebra is separable with
$\unmul \eq \comul$.
\end{defi}

A Frobenius algebra $(A,\mul,\eta,\comul,\eps)$ is separable iff
  \be
  \mul \circ (\ida \oti \mul) \circ (\ida \oti \zeta \oti \ida) \circ \Delta = \ida
  \label{eq:mulev3-3.1}
  \ee
for some morphism $\zeta \iN \Hom(\one,A)$ \Cite{Prop.\,3.1}{mulev3}.
In particular, for a connected Frobenius algebra, i.e.\ one satisfying 
$\Hom(\one,A) \eq \ko\,\eta$, separability and \deltasepy\ are equivalent.

Dually to Definition \ref{def:Ga-graded.alg}, a \emph{$\Ga$-graded coalgebra} in \C\
is a coalgebra $A \eq (\dt A,\comul,\eps)$ whose underlying object 
is $\Ga$-graded and for which the comultiplication $\comul$ is homogeneously even,
i.e.\ the morphisms $\comul^{i,j}_k \df (r^i \oti r^j) \,{\circ}\, \comul \cir e_k$ in
$\Hom(A_k,A_i\oti A_j)$ obey $\comul^{i,j}_k \eq 0$ for $k \nE i \pl j $.
As in the case of the unit morphism of an algebra, with the natural convention that 
$\one$ is $\Ga$-graded with $\one_j \eq \delta_{j,0}\, \one$, the counit morphism 
$\eps$ is automatically homogeneously even.
For a $\Ga$-graded coalgebra we abbreviate $\comul^{i,j}_{i+j} \,{=:}\, \comul^{i,j}$.
A \emph{$\Ga$-graded Frobenius algebra} in \C\ is a Frobenius algebra
$A \eq (\dt A,\mul,\eta,\comul,\eps)$ which is $\Ga$-graded both as an algebra
and as a coalgebra.

\begin{rem}
Graded Frobenius algebras in the monoidal category $\mathrm{comod}\text-H$ 
of comodules over a Hopf \ko-algebra $H$ have been considered in \cite{danN}.
If $H$ is the group algebra $\ko G$ of a group $G$, an algebra in 
$\mathrm{comod}\text-H$ is just a $G$-graded \ko-algebra. 
As shown in Corollary 4.2 of \cite{danN}, a (strongly) graded finite-dimensional
\ko-algebra $A$ is graded Frobenius if and only if its neutral component $A_0$ is 
Frobenius. For algebras $A$ in a generic \ko-linear monoidal category, being graded
Frobenius is still sufficient for $A_0$ being Frobenius, but the opposite implication
no longer holds.
\end{rem}

\begin{exa} \label{exa:twistedgroupalg2}
For $\Ga \,{\le}\, \picc$, the twisted group algebra $A_{(\Ga)}$ described in
Example \ref{exa:twistedgroupalg} can be endowed with the structure of a $\Ga$-graded
coalgebra by setting \Cite{Lemma\,3.10(iii)}{fuRs9} 
  \be
  \comul^{i,j} := |\Ga|\inv\coco(i,j)\inv\,\unlam^{i,j} ,
  \ee
where we write $\unlam^{i,j}$ for the basis of the one-dimensional vector space 
$\Hom(L_{i+j},L_i\oti L_j)$ that is covector dual to the basis $\lam_{i,j}$ of 
$\Hom(L_i\oti L_j,L_{i+j})$ used in the definition of the multiplication of $A_{(\Ga)}$.
This makes $A_{(\Ga)}$ into a \deltasep\ $\Ga$-graded Frobenius algebra.
\end{exa} 

\begin{exa} \label{exa:degenerate2}
We define the \emph{matrix representation}
$\comulM \eq \big( \comulb_{i,j} \big)_{\!i,j\in\Ga}$ of the comultiplication
of a $\Ga$-gra\-ded coalgebra with zero- or \onedim\ morphism spaces
$\Hom(A_{i+j},A_i\oti A_j)$ analogously
as the matrix representation of $\mu$ in Convention \ref{conv:mulmatrix}.
Similarly as in Example \ref{exa:degenerate1}, if $X \eq \one \oplus x$
is endowed with the comultiplication $\comul$ with matrix representation
  \be
  \comulM = \begin{pmatrix} 1 & 1 \\ 1 & 0 \end{pmatrix} ,
  \ee
then it is a \zz-graded counital coalgebra. However, this structure
does not fit together with the multiplication \eqref{eq:mulM.degenerate} into 
a Frobenius structure.
On the other hand, provided that $x$ is simple one does obtain a Frobenius algebra 
with multiplication \eqref{eq:mulM.degenerate} by endowing $X$ with a Frobenius form
that has the same matrix representation as $\mul$. However, this does not yield
a \zz-graded coalgebra structure on $X$.
\end{exa}

\begin{lem} \label{def:lem:twist-Frobalg}
Let \C\ be a monoidal category and $A \eq (\dt A,\mul,\eta,\comul,\eps)$ be a
$\Ga$-graded Frobenius algebra in \C.
Then $A\ococ \df (\dt A,\mulco,\eta,\comulco,\eps)$ with cocycle-twisted
multiplication \eqref{eq:def:wmulijk} and co\-cy\-cle-twisted comultiplication
$\comul^{\!\!(\coc)\,i,j}_{} \df \coc(i,j)\inv \comul^{i,j}_{}$,
and with the same unit and counit as $A$, is again a Frobenius algebra in \C.
Further, if $A$ is separable, then so is $A\ococ$.
\end{lem}

\begin{proof}
That $A\ococ$ is an algebra is explained in Remark \ref{rem:twist-alg};
that it is a coalgebra follows in the same way. Concerning the Frobenius property,
by direct calculation we have (taking again for simplicity \C\ to be strict)
  \be
  \begin{aligned}
  (\ida \oti \mulco) \cir (\comulco \oti \ida)
  = \sum_{i,j,k\in\Ga}\! & \coc(i{+}k,-k)\inv \coc(-k,j)\, (e_{i+k} \oti e_{j-k})
  \\[-6pt]
  & \cir (\id_{A_{i+k}} \oti \mul_{-k,j})
  \cir (\comul^{i+k,-k} \oti \id_{A_j}) \cir (r^i \oti r^j)
  \end{aligned}
  \label{eq:comp:wFrob1}
  \ee
and
  \be
  \comulco \cir \mulco 
  = \sum_{i,j,k\in\Ga}\! \coc(i{+}k,j{-}k)\inv \coc(i,j)\; (e_{i+k} \oti e_{j-k})
  \circ \comul^{i+k,j-k} \cir \mul_{i,j} \circ (r^i \oti r^j) \,.
  \label{eq:comp:wFrob2}
  \ee
Since $\coc$ is a two-cocycle, the products of $\coc$-factors in these two
expressions are equal, and hence equality of \eqref{eq:comp:wFrob1} and
\eqref{eq:comp:wFrob2} is implied by the corresponding equality for the untwisted 
operations, i.e.\ by the Frobenius property of $A$.
The second Frobenius identity follows analogously.
The statement about separability follows directly from the fact that
$\mulco_{i,j} \cir \comul^{\!\!(\coc)\,i,j}_{} \eq \mul_{i,j} \cir \comul^{\!i,j}_{}$
together with the characterization \eqref{eq:mulev3-3.1} of separability.
\end{proof}

\begin{rem}
Like for \ko-algebras, the property of admitting a Frobenius structure is not Morita
invariant. Rather, the Morita class of a Frobenius algebra in a monoidal category
may contain quasi-Frobenius algebras that are not Frobenius \cite{shimi27}.
In contrast, admitting a \emph{symmetric} Frobenius structure \emph{is}
a Morita invariant property. The notion of symmetric Frobenius algebra, which is
defined in any pivotal monoidal category, is important in applications like those
in \cite{muge8,fuRs4,henri}. This notion does not play a role for the considerations 
in the present paper, but it is worth noting that the
full subcategory generated by the homogeneous components of a twisted group
algebra, being equivalent to $\VectGe$, admits a pivotal structure.
\end{rem}

%%%%%%%%%%%%%%%%%%%%%%%%%%%%%%%%%%%%%%%%%%%%%%%%%%%%%%%%%%%%%%%%%%%%%%%%

\subsection{Tensor product of bimodules} \label{sec:prel-3}

The standard tensor product of bimodules over an algebra in the category
of \ko-vector spaces readily generalizes to other \ko-linear categories:

\begin{defi} \label{def:ota-coequ}
Let $A$ be an algebra in a strict monoidal category \C, and let $m \eq (\dtm,\ohr)$
and $n \eq (\dtn,\rho)$ be a right and a left $A$-module in \C, respectively.
If the coequalizer
  \be
  m \ota n := \mathrm{coequ} \big(\!
  \begin{tikzcd}[column sep=3.1em]
  \dtm \oti A \oti \dtn \ar[yshift=3pt]{r}{\ohr\oti\id_\dtn\,}
  \ar[yshift=-2pt]{r}[swap]{\id_\dtm\oti\rho\,}
  & \dtm \oti \dtn
  \end{tikzcd}
  \!\big) \,,
  \label{eq:ota-as-coequ}
  \ee
exists, then it is called the \emph{tensor product over} $A$ of $m$ and $n$.
\end{defi}

Note that $\id_\dtm \oti \eta \oti \id_\dtn$ is a common section of
$\ohr\oti\id_\dtn$ and $\id_\dtm\oti\rho$, i.e.\ the coequalizer 
\eqref{eq:ota-as-coequ} is reflexive.
Also, here and analogously for the statements below, strictness of \C\ is assumed only
for notational convenience; the only modification in the non-strict case is that
one of the two morphisms that are coequalized gets composed with an associator.

Assume now that \C\ admits all coequalizers of the form \eqref{eq:ota-as-coequ} and
that the tensor product functor $\oti$ of \C\ preserves such coequalizers. Then 
for any two $A$-bimodules $m$ and $n$ 
the object $m \ota n$ admits a natural structure of $A$-bimodule as well. Hence the
tensor product \eqref{eq:ota-as-coequ} endows the category $A\bimod$ of $A$-bimodules
with the structure of a monoidal category, with $A$ as monoidal unit and with 
associator and unitors inherited from those of \C; compare \Cite{Sect.\,4.1}{gaJoy}
for a detailed exposition in a bicategorical context, and e.g.\
\Cite{Exc.\,7.8.22,7.8.28}{EGno} for the case that \C\ is a multitensor category.

If 
\C\ is abelian and
the algebra $A$ is separable Frobenius, then the tensor product $\otA$ can be 
realized as an image. Denote by
  \be
  \coe \Colon \dtm \oti \dtn \rarr~ m \ota n
  \label{eq:r:mn2motan}
  \ee
the epimorphism given by the definition of the coequalizer \eqref{eq:ota-as-coequ}.
Then, as a moderate generalization of Lemma\,1.21 and Corollary 1.22 of 
\cite{kios}, we have

\begin{lem} \label{lem:mAn=Im(p)}
Let $A$ be an algebra in
an abelian
strict monoidal category \C, and let 
$m \eq (\dtm,\ohr)$ be a right $A$-module and $n \eq (\dtn,\rho)$ a left $A$-module.
If $\pp \,{\equiv}\, \pp_{m,n} \iN \End(\dtm\oti\dtn)$ is an idempotent and obeys 
$\coe \cir \pp \eq \coe$
as well as
  \be
  \pp \circ (\ohr\oti\id_{\dtn}) = \pp \circ (\id_{\dtm}\oti\rho) \,,
  \label{eq:p1prop}
  \ee
then we have
  \be
  m \ota n \cong \Im(\pp) \,,
  \ee
and there is a direct sum decomposition $\dtm\oti\dtn \Cong  m \ota n \,{\oplus}\, \Ker(\pp)$.
Moreover, if in addition $A$ is separable Frobenius with separability morphism
$\unmul$, then the endomorphism
  \be
  p \equiv p_{m,n}
  := (\ohr\oti\rho) \circ (\id_{\dtm} \oti (\unmul\cir\eta) \oti \id_{\dtn})
  \label{eq:def:bimod-p}
  \ee
satisfies these three properties.
\end{lem}

\begin{proof}
By the universal property of $m\ota n$ as a coequalizer, validity of
\eqref{eq:p1prop} implies that the kernel $\ker(\coe)$ is a subobject of 
$\ker(\pp)$, while $\coe \cir p \eq \coe$ implies that $\ker(\pp)$ is a subobject of
$\ker(\coe)$. Thus $\ker(\pp) \eq \ker(\coe)$ and hence, using that $\coe$ is epi, 
$\Im(\pp) \Cong \Im(\coe) \eq m \ota n$.
 \\
Further, using that $A$ is Frobenius it follows that the morphism $p$ given by
\eqref{eq:def:bimod-p} satisfies \eqref{eq:p1prop}, and using that $A$ is separable
Frobenius it follows that $p$ is an idempotent.
Finally, by the definition of $m\ota n$ as a coequalizer we have
  \be
  \coe \circ p = \coe \circ \big[ \id_{\dtm} \oti \big(\rho \cir (\ida\oti\rho) \cir
  ((\unmul\cir\eta)\oti \id_{\dtn})\big) \big]
  \ee
which when combined with the separability and the representation property of $\rho$
simplifies to $\coe \cir p \eq \coe$.
\end{proof}

\begin{rem}
The requirement that the category \C\ (which by our general assumptions is additive)
is abelian can be relaxed. It suffices that the kernel and image of the
epimorphism \eqref{eq:r:mn2motan} and of the idempotent $\pi$ exist and that
the idempotent $\id_{\dtm\otimes\dtn}\,{-}\,\pi$ is again split.
\end{rem}

Thus in short, if $A$ is separable Frobenius, then we have
$m \ota n \Cong \Im(p)$ with $p \,{\equiv}\, p_{m,n}$ as given in 
\eqref{eq:def:bimod-p}. Moreover, the tensor product of morphisms is then given by
  \be
  f \ota g = p_{m',n'} \circ (f \oti g) = (f \oti g)  \circ p_{m,n}
  \label{eq:ota=oti.p}
  \ee
for $f\colon m \Rarr~ m'$ and $g\colon n \Rarr~ n'$ morphisms of right and left
$A$-modules, respectively.

%%%%%%%%%%%%%%%%%%%%%%%%%%%%%%%%%%%%%%%%%%%%%%%%%%%%%%%%%%%%%%%%%%%%%%%%

\subsection{\Graco\ algebras} \label{sec:prel-4}

Let now \C\ in addition be braided, with braiding $\br$. Given an algebra
$A \eq (\dt A,\mul,\eta)$ in \C, the \emph{opposite algebra} $\Ao$, defined as
$\Ao \,{:=}\, (\dt A,\mul \cir \br_{\dt A,\dt A},\eta)$,
is again a unital associative algebra in \C.

\begin{lem} \label{lem:Wm}
Let $A$ be an algebra in a braided monoidal category.
The categories $A\Mod$ of left $A$-modules and $\moD\Ao$ of right $\Ao$-modules 
are isomorphic. Specifically, if $n \eq (\dtn,\rho)$ is a left $A$-module, then 
$(\dtn,\rho \cir \br_{\dtn,A})$ is a right $\Ao$-module, and if $(\dtm,\ohr)$ is
a right $\Ao$-module, then $(\dtm,\ohr\cir\br_{A,\dtm})$ is a left $A$-module.
\end{lem}

\begin{proof}
That, given a left $A$-module $(\dtn,\rho)$, the morphism $\rho \cir \br_{\dtn,A}$
defines a right $\Ao$-action on $\dtm$ follows directly with the help of the
naturality of the braiding. 
In terms of the standard graphical string calculus
for monoidal categories, the compatibility with the multiplication on $\Ao$ is 
seen as follows:
\def\locC    {0.3}   %  radius circle
\def\locH    {3.4}   %  height straight line
\def\locR    {2.8}   %  height upper rep symbol
\def\locS    {1.5}   %  height lower rep symbol
\def\locW    {0.5}   %  hor distance
  \be
  \raisebox{-4.7em}{
  \begin{tikzpicture}[scale=\locscale]
  \draw[line width=\widthObj,\colorAlg]
       (2*\locW,0) node[below=1pt] {$A$} -- +(0,0.45*\locR)
       [out=90,in=270] to (-0.8*\locW,0.85*\locR) [out=90,in=220] to (0,\locR)
       (\locW,0) node[below=1pt] {$A$} -- +(0,0.3*\locS)
       [out=90,in=270] to (-0.7*\locW,0.7*\locS) [out=90,in=220] to (0,\locS) ;
  \draw[line width=\widthBcg,white] (0,0) -- +(0,\locH) ;
  \draw[line width=\widthObj,\colorMod]
       (0,0) node[below=1pt] {$\dtn$} -- +(0,\locH) node[above=-1pt] {$\dtn$} ;
  \drawRhoLabel {0}{\locR} \rho ;
  \drawRhoLabel {0}{\locS} \rho ;
  \end{tikzpicture}
  }
  = ~
  \raisebox{-4.7em}{
  \begin{tikzpicture}[scale=\locscale]
  \draw[line width=\widthObj,\colorAlg]
       (2*\locW,0) node[below=1pt] {$A$} -- +(0,0.15*\locR)
       [out=90,in=270] to (-3*\locC,\locR-3*\locC) arc (180:90:\locC)
       [out=90,in=220] to (0,\locR)
       (\locW,0) node[below=1pt] {$A$} -- +(0,0.05*\locR)
       [out=90,in=270] to (-3*\locC,\locR-5.5*\locC) ;
  \draw[line width=\widthBcg,white]
       (0,0) -- +(0,\locH)
       (-3*\locC,\locR-5.5*\locC) [out=90,in=270]to  (-\locC,\locR-3*\locC) ;
  \draw[line width=\widthObj,\colorAlg]
       (-3*\locC,\locR-5.5*\locC) [out=90,in=270] to (-\locC,\locR-3*\locC)
       arc (0:90:\locC) ;
  \draw[line width=\widthObj,\colorMod]
       (0,0) node[below=1pt] {$\dtn$} -- +(0,\locH) node[above=-1pt] {$\dtn$} ;
  \drawRhoLabel {0}{\locR} \rho ;
  \drawMulLabelL {-2*\locC}{\locR-2*\locC} \mul ;
  \end{tikzpicture}
  }
  = ~
  \raisebox{-4.7em}{
  \begin{tikzpicture}[scale=\locscale]
  \draw[line width=\widthObj,\colorAlg]
       (\locW+2*\locC,0) node[below=1pt] {$A$} -- +(0,0.15*\locR)
       [out=90,in=270] to (\locW,0.4*\locR) arc (180:90:\locC)
       [out=90,in=270] to (-0.8*\locW,0.85*\locR) [out=90,in=220] to (0,\locR) ;
  \draw[line width=\widthBcg,white]
       (0,0) -- +(0,\locH)
       (\locW,0) -- +(0,0.15*\locR) [out=90,in=270]
       to (\locW+2*\locC,0.4*\locR) arc (0:90:\locC) ;
  \draw[line width=\widthObj,\colorAlg]
       (\locW,0) node[below=1pt] {$A$} -- +(0,0.15*\locR) [out=90,in=270]
       to (\locW+2*\locC,0.4*\locR) arc (0:90:\locC) ;
  \draw[line width=\widthObj,\colorMod]
       (0,0) node[below=1pt] {$\dtn$} -- +(0,\locH) node[above=-1pt] {$\dtn$} ;
  \drawRhoLabel {0}{\locR} \rho ;
  \drawMulLabelR {\locW+\locC}{0.4*\locR+\locC} \mul ;
  \end{tikzpicture}
  }
  \ee
(Our conventions for the string calculus can be easily read off from these
pictures.)
Moreover, module morphisms with respect to $\rho$ are also module morphisms with respect 
to $\rho \cir \br_{\dtm,A}$, as follows again from naturality of the braiding. One
therefore deals with a functor from $A\Mod$ to $\moD\Ao$. By invertibility of the
braiding, this functor is an isomorphism.
\end{proof}

\begin{rem}
By iteration, for any $\ell\iN\zet$,
$A^{(\ell)} \df (\dt A,\mul \cir (\br_{\dt A,\dt A})^\ell_{\phantom:},\eta)$
is an algebra structure on the object $\dt A$. Via appropriate even and odd powers,
respectively, of the braiding, all categories $A^{(2r)}\Mod$ and $\moD A^{(2r+1)}$ 
for $r\iN\zet$ are isomorphic.
Also, if $A \eq (\dt A,\mul,\eta,\comul,\eps)$ is a Frobenius algebra, then
  \be   
  A^{(\ell)} := (\dt A,\mul \cir (\br_{\dt A,\dt A})^\ell_{\phantom|},\eta,
  (\br_{\dt A,\dt A})^{-\ell}_{\phantom|}{\circ}\, \comul,\eps)
  \ee
is again a Frobenius algebra structure on the object $\dt A$.
Moreover, if $A$ is separable, then so is $A^{(\ell)}$, and if 
if $A$ is \deltasep, then so is $A^{(\ell)}$.
\end{rem}

A \emph{commutative} algebra in a braided monoidal category \C\ is an algebra $A$
for which the identity map is an isomorphism between $A$ and
its opposite algebra $\Ao$ (and thus also all algebras $A^{(n)}$).
In terms of the components \eqref{eq:mul-ijk} of the multiplication, this means
$\mul_{j,i} \cir \br_{A_i,A_j} \eq 
\mul_{i,j}$.
For graded algebras we give analogously

\begin{defi} \label{def:graco}
Let \C\ be a braided monoidal category. A $\Ga$-graded algebra
$A \eq (\dt A,\mul,\eta)$ in \C\ is called \emph{\graco} iff
  \be
  \mul_{j,i} \circ \br_{A_i,A_j} = \coc(i,j)\, \mul_{i,j}
  \label{eq:def:graco}
  \ee
for $i,j\iN\Ga$ for a normalized two-cocycle $\coc$ on $\Ga$.
\end{defi}

A commutative algebra is thus the same as a \graco\ algebra for $\coc$ the trivial
two-cocycle. Note that in our terminology `\graco' we neither specify the relevant 
two-cocycle nor the relevant braiding explicitly; this is legitimate because in the
situations we are interested in, we deal with a given fixed cocycle $\coc$ and a 
given fixed braiding $\br$.

\begin{exa} \label{exa:twistedgroupalg3}
Consider again the twisted group algebras from Example \ref{exa:twistedgroupalg},
with underlying objects $A_{(\Ga)} \eq \bigoplus_{i\in\Ga\!}L_i$, and let now \C\ be
braided (so that $\picc$, and hence $\Ga$, is necessarily abelian).
The full subcategory $\C_{(\Ga)}$ of \C\ whose objects are direct sums of the $L_i$
with $i\iN\Ga$ is equivalent to the semisimple abelian category $\VectGe$ of 
finite-dimensional $\Ga$-graded \ko-vec\-tor spaces with homogeneously even linear 
maps. The inequivalent associators and braidings of a categorical group, and thus also
those for the Picard category $\Picc$ and its full monoidal subcategory $\C_{(\Ga)}$,
are classified by the third abelian cohomology \Cite{Sect.\,3}{joSt6}. (An
abelian three-cocycle $(\Psi,\Coco) \iN Z^3_{\rm ab}(G,\kox)$ on a group $G$ 
consists of an ordinary group three-cocycle $\Psi \iN Z^3(G,\kox)$ and a two-cochain 
$\Coco \iN C^2(G,\kox)$ that satisfies compatibility relations with $\Psi$, see
\cite{joSt6}.)
An abelian three-cocycle $(\Psi,\Coco) \iN Z^3_{\rm ab}(\picc,\kox)$ encodes
precisely the restriction to $\Picc$ of the associator and braiding of \C. In terms 
of bases $\lam_{i,j}$ of $\Hom(L_i\oti L_j,L_{i+j})$ (cf.\ Example 
\ref{exa:twistedgroupalg}), the braiding on $\Picc$ is given by
$\lam_{j,i} \cir \br_{L_i,L_j} \eq \Coco(j,i)\inv \lam_{i,j}$, while the 
associator sends $\lam_{i+j,k} \cir (\lam_{i,j} \oti \id_{L_k})$ to
$\Psi(i,j,k)\inv (\lam_{i,j+k} \cir (\id_{L_i} \oti \lam_{j,k}))$
\Cite{App.\,C.2}{fuRs9}. Two abelian three-cocycles yield equivalent braided 
structures iff they are cohomologous. Further \Cite{Prop.\,3.14(i)}{fuRs9}, every 
simple special Frobenius algebra in $(\C_{(\Ga)};\psi,\coco)$, with
$(\psi,\coco) \eq (\Psi,\Coco)|_{\Ga}^{}$, is isomorphic to one of the twisted
group algebras $A_{(\Ga)}$ with multiplication \eqref{eq:mul=coco.lam}.
\end{exa}

\begin{lem} \label{lem:graco-iff-bichar}
The twisted group algebra $A_{(\Ga)}$ in $(\C_{(\Ga)};\psi,\coco)$ is \graco\ iff 
the cochain $\coco$ is a bicharacter on $\Ga$.
\end{lem}

\begin{proof}
Combining the formula \eqref{eq:mul=coco.lam} for the multiplication 
with the defining relation \eqref{eq:def:graco} for a \graco\ algebra 
and with the expression of the braiding in terms of $\coco$,
one finds that $\coco \eq \coc\inv$. 
In particular, the two-cochain $\coco$ must in fact be a two-cocycle, and as a
consequence $\psi \eq \mathrm d\coc \eq 1$. Further, by definition of $Z^3_{\rm ab}$
one has $(1,\coco) \iN Z^3_{\rm ab}(\Ga,\kox)$ iff $\coco$ is a bicharacter,
i.e.\ satisfies $\coco(i{+}j,k) \eq \coco(i,k)\,\coco(j,k)$
and $\coco(i,j{+}k) \eq \coco(i,j)\,\coco(i,k)$ for $i,j,k\iN\Ga$.
\end{proof}

\begin{exa} \label{exa:zz-2}
Consider the \zz-graded Frobenius algebra $A \eq \one \,{\oplus}\, L$ with $L$ an 
order-2 invertible object satisfying $\FxL \df \FO L \eq 1$, as described in Example 
\ref{exa:zz}, in a monoidal category \C\ with Picard group \zz. If \C\ is braided,
then the category $\Picc$ is braided as well.
The elements of the abelian cohomology group $H^3_{\rm ab}(\zz,\kox)$ can be labeled
by the possible values of $\FxL$ and of the multiple of $\id_{L\otimes L}$ that 
constitutes the self-braiding $\br_{L,L}$. Compatibility of associator and 
braiding requires that the latter number squares to $\FxL$, so that in the situation 
at hand we have $\br_{L,L} \eq {\pm}\, \id_{L\otimes L}$.
For $\br_{L,L} \eq \id_{L\otimes L}$, $\Picc$ is braided equivalent to the category 
of finite-di\-mensional \zz-graded vector spaces and $A \eq \one \,{\oplus}\, L$ is a 
commutative algebra in \C, while for $\br_{L,L} \eq {-}\,\id_{L\otimes L}$, $\Picc$ 
is braided equivalent to the category of finite-dimensional supervector spaces 
and $A$ is a supercommutative algebra in \C.
In the particular case of $\C \eq \C(\mathrm A_1^{(1)},k)$, as defined in
Example \ref{exa:zz}, one finds $\br_{L,L} \eq \mathrm i^k\,\id_{L\otimes L}$ 
\Cite{Exc.\,8.18.6}{EGno}, hence for $k \iN 4\zet$ the object $\one \,{\oplus}\, U_k$
admits a natural structure of a commutative Frobenius algebra, while for
$k\iN4\zet{+}2$ it admits a natural structure of a supercommutative Frobenius algebra.
\end{exa}

\begin{exa} \label{exa:degenerate3}
The \zz-graded algebra $X \eq \one \oplus x$ from Example \ref{exa:degenerate1}
is \graco\ with respect to any braiding on \C.
\end{exa}

%%%%%%%%%%%%%%%%%%%%%%%%%%%%%%%%%%%%%%%%%%%%%%%%%%%%%%%%%%%%%%%%%%%%%%%%

\subsection{Monoidal monads} \label{sec:prel-5}

Pertinent aspects of the representation theory of an algebra $A$ in a monoidal 
category \C\ can conveniently be studied in terms of the corresponding monad on \C,
whose underlying endofunctor is $A\oti{-}$. We will use this approach in Section
\ref{sec:gil}.

Recall that a \emph{monad} on a category \C\ is an algebra $\TT \eq (T,\mulT,\etaT)$
in the monoidal category of endofunctors of \C. Thus $\mulT$ and $\etaT$ are natural 
transformations $\mulT\colon T^2 \RArr~ T$ and $\etaT\colon \Id_\C \RArr~ T$ satisfying
  \be
  \mulT \circ \mulT\,T = \mulT \circ T\mulT \qquand
  \mulT \circ \etaT\,T = \Id_\C  = \mulT \circ T\etaT ,
  \label{eq:mmT=mTm}
  \ee
or, in components, $\mulT_x \cir \mulT_{Tx} \eq \mulT_x \circ T\mulT_x$ and 
$\mulT_x \circ \etaT_{Tx} \eq \id_{Tx} \eq \mulT_x \cir T\etaT_x$ for $x\iN\C$.
Naturality of $\mulT$ and $\etaT$ means that $\mulT_{Ty} \cir TTf \eq Tf \cir \mulT_x$
and that $Tf \cir \etaT_x \eq \etaT_y \cir f$ for $f\colon x \Rarr~ y$, respectively.

\begin{defi}
A \emph{module} over a monad $\TT$ in \C\ is a pair $(x,\rho_x)$ consisting of
an object $x \iN \C$ and a morphism $\rho_x\colon Tx \Rarr~ x$ satisfying
$\rho_x \circ \mulT_x \eq \rho_x \circ T\rho_x$ and
$\rho_x \circ \etaT_x \eq \id_x$.
A \emph{morphism of $\TT$-modules} from $(x,\rho_x)$ to $(y,\rho_y)$ is a morphism
$f\colon x \Rarr~ y$ in \C\ satisfying $f \cir \rho_x \eq \rho_y \cir Tf$.
\end{defi}
	
In part of the literature a $\TT$-module is instead called a $\TT$-\emph{algebra}.
In the present context this terminology would be inconvenient.

\begin{defi} \label{def:monoidalmonad}
A \emph{monoidal structure} on a monad $\TT \eq (T,\mulT,\etaT)$ on a monoidal 
category \C\ is a family $\tZ^{} \,{\equiv}\, \tZ^\TT$ of morphisms
  \be
  \TZ xy \Colon Tx \oti Ty \rarr~ T(x \otI y)
  \ee
that are natural in $x,y \iN \C$ and make $(T,\tZ,\etaT_\one)$ into a 
lax monoidal functor and $\mulT$ and $\etaT$ into monoidal natural transformations. 
\end{defi}

Put differently, a monoidal monad is a monad in the 2-category of monoidal categories,
lax monoidal functors and monoidal natural transformations.
Explicitly, the conditions that the natural transformations $\mulT$ and $\etaT$ 
are monoidal are given by the commutativity of the diagrams
  \be
  \begin{tikzcd}[column sep=3.1em, row sep=2.4em]
  T^2x \oti T^2y \ar{r}{\TZ{Tx}{Ty}} \ar{dr}[swap]{\mulT_x \oti \mulT_y}
  & T(Tx \oti Ty) \ar{r}{T(\TZ xy)}
  & T^2(x \otI y) \ar{d}{\mulT_{x\otimes y}}
  \\
  ~ & Tx \oti Ty \ar{r}{\TZ xy} & T(x\otI y)
  \end{tikzcd}
  \label{eq:mulT.montrafo}
  \ee
and
  \be
  \begin{tikzcd}[column sep=3.5em, row sep=2.2em]
  x \oti y \ar{r}{\etaT_x \oti \etaT_y\,}
  \ar{dr}[swap,xshift=4pt,yshift=4pt]{\etaT_{x\otimes y}}
  & Tx \oti Ty \ar{d}{\TZ xy}
  \\
  & T(x\otI y)
  \end{tikzcd}
  \ee
respectively.

Dually to the notion of a monad, a \emph{comonad} on a category \C\ is a coalgebra 
$\TT \eq (T,\comulT,\epsT)$ in the category of endofunctors of \C, i.e.\ the natural 
transformations $\comulT\colon T\RArr{~} T^2$ and $\epsT\colon T \RArr{~} \Id_\C$
satisfy $\comulT T \cir \comulT \eq T \comulT \cir \comulT$ and
$\epsT T \cir \comulT \eq \Id_\C \eq T \epsT \cir \comulT$.

\begin{defi} \cite{stre8,doPe9} \label{def:frobmonad}
(i) A \emph{Frobenius monad} is a monad $(T,\mulT,\etaT)$ together with a comonad
structure $(T,\comulT,\epsT)$ such that the equalities 
  \be
  (\mulT \oti \id_T)  \circ (\id_T \oti \comulT) = \comulT \circ \mulT
  = (\id_T \oti \mulT) \circ (\comulT \oti \id_T)
  \label{eq:frobmonad}
  \ee
hold.
 \\[2pt]
(ii) A \emph{commutative Frobenius monad} is a Frobenius monad $\TT$ together with a
natural isomorphism $\tau \colon T^2 \RArr~ T^2$ satisfying the hexagon identity
$T\,\tau \cir \tau\,T \cir T\,\tau \eq \tau\,T \circ T\,\tau \circ \tau\,T$
and being compatible with the monad and comonad structures according to
$\mulT \cir \tau \eq \mulT$, $\tau \cir \etaT\,T \eq T\,\etaT$,
$\tau \cir \comulT \eq \comulT$ and $\epsT T \cir \tau \eq T\,\epsT$, as well as
$\tau \cir \mulT T \eq T \mulT \cir \tau\,T \cir T\,\tau$ and
$\comulT T \cir \tau \eq T\,\tau \cir \tau\,T \cir T \comulT$.
 \\[2pt]
(iii) A \emph{separable} monad is a monad whose multiplication has a one-sided 
inverse $\unmul^\TT$, so that $\mulT \cir \unmul^\TT \eq \id_T$
 \\[2pt]
(iv) A \emph{\deltasep} monad is a monad-comonad $(T,\mulT,\etaT,\comulT,\epsT)$
which is separable with $\unmul^\TT \eq \comulT$, so that $\mulT \cir \comulT \eq \id_T$.
\end{defi}

\begin{rem}
Just like for Frobenius algebras, there are several equivalent definitions of 
Frobenius monad \Cite{Thm.\,1.6}{stre8}. They basically boil down to the statement
that the Frobenius pairing $\epsT \cir \mulT$ and Frobenius copairing 
$\comulT \cir \etaT$ furnish the counit and unit of a two-sided self-adjunction 
$T \,{\dashv}\, T \,{\dashv}\, T$. To appreciate this, recall that every 
adjunction $\adj$ of functors $F\colon \C \,{\leftrightarrows}\, \D \colon G$ 
($G$ right adjoint to $F$) gives rise to a monad $(T,\mulT,\etaT)$ on \C,
with $T \eq GF$, unit $\etaT$ being the unit $\Eta$ of the adjunction,
and multiplication given by whiskering the counit $\Eps$ of the adjunction,
$\mulT \eq \id_G\,\Eps\,\id_F$ \Cite{Lemma\,5.1.3}{RIeh2}.
Likewise, it gives rise to a comonad $(S,\comulS,\epsS)$ on \D, with $S \eq FG$,
$\epsS \eq \Eps$ and $\comulS \eq \id_F\,\Eta\,\id_G$.
If the adjunction $\adj$ is two-sided, then there is accordingly both a monad
and a comonad structure with underlying endofunctor $T$ (and likewise for $S$);
these fit together to form a Frobenius monad (this is easily seen when using
bicategorical string diagrams, see e.g.\ \Cite{Sect.\,1.2}{muge8}).
\end{rem}

To any algebra $A \eq (\dt A,\mul,\eta)$ in a monoidal category \C\ there is
naturally associated a monad $\AA \eq (T_A,\mulA,\etaA)$ on \C, with
$T_A \eq {-}\oti \dt A$, and with multiplication and unit given by
natural transformations $\mulA$ and $\etaA$ with components
  \be
  \mulA_c = (\id_c \oti \mul) \cir \ass_{c,A,A} \qquand
  \etaA_c \eq (\id_c \oti \eta) \cir \Unr c\inv
  \label{eq:def:AA}
  \ee
for $c\iN\C$. Indeed, the defining properties of $\mulA$ and $\etaA$ are 
equivalent to the associativity and unit properties of $\mul$ and $\eta$.
Likewise, to any coalgebra $A \eq (\dt A,\comul,\eps)$ in \C\ there is
associated a comonad $(T_A,\comulA,\epsA)$ on \C, with underlying endofunctor $T_A$ 
and with natural transformations having components
  \be
  \comulA_c = \ass_{c,A,A}\inv \cir (\id_c \oti \comul) \qquand
  \epsA_c = \Unr c \cir (\id_c \oti \eps) \,.
  \label{eq:def:coAA}
  \ee
By direct verification of the required equalities one checks:

\begin{lem} \label{lem:A-AA-Frobenius}
Let $A$ be a Frobenius algebra in a monoidal category \C. The induced monad and
comonad structures on the endofunctor $T_A \eq {-}\oti \dt A$ fit together to the 
structure of a Frobenius monad $\AA \eq (T_A,\mulA,\etaA,\comulA,\epsA)$ on $T_A$.
If $A$ is a separable (\deltasep) algebra, then $\AA$ is a separable
(\deltasep) monad.
If \C\ is braided and $A$ is commutative, then $\AA$ is a commutative Frobenius 
monad, with $\tau_c \eq \id_c \oti \br_{\dt A,\dt A}$.
\end{lem}

%%%%%%%%%%%%%%%%%%%%%%%%%%%%%%%%%%%%%%%%%%%%%%%%%%%%%%%%%%%%%%%%%%%%%%%%

\subsection{Graded-enriched categories} \label{gr-enriched}

For $\Ga$ a finite abelian group, the category \VectG\ has
as objects $\Ga$-graded \ko-vector spaces $V \eq \bigoplus_{i\in\Ga\!} V_i$
and as morphisms all \ko-linear maps. The category $\VectGe$ is its wide subcategory 
whose morphisms are the grade-preserving \ko-linear maps.
A symmetric monoidal structure on \VectG\ that inherits the associativity and unit 
constraints from the monoidal category $(\Vect,\otimes)$ of \ko-vector spaces 
(and which we thus still denote by $\otimes$) is obtained by setting
$(V \oti V')_i \,{:=}\, \bigoplus_{j\in\Ga} V^{}_j \oti V'_{i-j}$ for objects and
$(f \oti f')(v \oti v') \,{:=}\, f(v) \otimes f'(v')$
for morphisms of \VectG. The subcategory $\VectGe$ inherits this symmetric 
monoidal structure.

As is well known, for $\Ga \eq \zz$ the monoidal structure can be twisted by the
non-trivial \zz-bicharacter (see e.g.\ \cite{brEl}), thereby yielding instead a 
\emph{supermonoidal} structure. A natural generalization of this construction to 
arbitrary finite abelian groups has been presented in \cite{puty,brDa}. 
In the sequel we follow the exposition in \cite{scVa}.

We denote the grade of a homogeneous vector $v \iN V_i$ by $\gr v \eq i$
and the grade of a homogeneous linear map $f$ by $\gr f$.

\begin{defi} \label{def:oticoc}
Let $\coc$ be a $\Ga$-bicharacter. The \emph{$\coc$-\twi\ tensor product} $\otIcoc$ 
on \VectG\ is defined on objects as $V \oticoc V' \,{:=}\, V \oti V'$
and on morphisms as the quadrilinear extension of the prescription
  \be
  (f \oticoc f')(v \oticoc v') := \coc(\gr{f'},\gr v)\, f(v) \oticoc f'(v')
  \label{eq:f*f'(v*v')}
  \ee
for homogeneous morphisms $f\colon V\Rarr~ W$, $f'\colon V'\Rarr~ W'$ and homogeneous
vectors $v\iN V$, $v'\iN V'$.
\end{defi}

Note that, as we always work with some fixed group $\Ga$, the dependence on
$\Ga$ is suppressed in the notation and terminology.
The $\coc$-\twi\ tensor product is \emph{not} a monoidal structure on \VectG. 
However, it enjoys properties analogous to those of a monoidal structure, except for
the following:

\begin{lem} \label{lem:GILvectG}
$($\emph{$\coc$-Twisted Interchange Law.}$)$ 
We have
  \be
  (f' \oticoc g') \circ (f \oticoc g)
  = \coc(\gr{g'},\gr f)\, (f' \cir f) \oticoc (g' \cir g)
  \label{eq:interchangelaw4vectG}
  \ee
for any quadruple of homogeneous morphisms $f,f',g,g'$ in \VectG.
\end{lem}

\begin{proof}
This follows immediately by combining the character property of $\coc$ in both of 
its arguments with the identities $\gr{f(v)} \eq \gr f\pl \gr v$ and 
$\gr{g'\cir g} \eq \gr{g'}\pl\gr g$.
\end{proof}

\begin{defi}\label{def:cocgracat}
  \Enumerate
  \item  \label{item:cocgracat}
A \emph{$\Ga$-graded \ko-linear category} is a category enriched over the 
symmetric monoidal category $(\VectGe,\otimes)$. 
  \item
A \emph{$\Ga$-graded \ko-linear functor} is a functor enriched over $(\VectGe,\otimes)$. 
  \item
We denote by \CatG\ the category having small $\Ga$-graded \ko-linear categories 
as objects and $\Ga$-graded \ko-linear functors as morphisms.
  \item
The \emph{$\coc$-\twi\ Cartesian product} $\C \boticoc \D$
of two $\Ga$-graded \ko-linear categories\,% 
 \footnote{~%
	In \cite{scVa} the symbol $\boxtimes^{\coc}$ is used for this product.}
is defined on objects as the Cartesian product, i.e.
  \be
  \obj(\C \boticoc \D) := \obj(\C) \Times \obj(\D)
  \ee
for $\C,\D\iN\CatG$, and on morphism spaces as the tensor product of vector spaces 
with composition of morphisms defined as in \eqref{eq:interchangelaw4vectG}. 
  \end{enumerate}
\end{defi}

We indicate the prescription for the composition in $\C \boticoc \D$ by writing
  \be
  \Hom_{\C \botIcoc \D}((c,d),(c',d')) = \HomC(c,c') \oticoc \HomD(d,d') \,.
  \ee
By construction, $\C \boticoc \D$ is endowed with the structure of a $\Ga$-graded 
\ko-linear category and $\botIcoc$ becomes a monoidal structure on the category \CatG.

\begin{defi}\label{def:cocgramoncat} \Cite{Def.\,2.4}{scVa}
A \emph{\cgmc} is a $\Ga$-graded \ko-linear category \C\
together with a unit object and a $\Ga$-graded \ko-linear functor
  \be
  \oticoc \Colon \C \boticoc \C \rarr~ \C
  \ee
that satisfy the
same unitality and associativity axioms as in the case of a strict monoidal category.
\end{defi}

In particular, the $\coc$-\twi\ tensor product $\otIcoc$ introduced in Definition 
\ref{def:oticoc} furnishes a $\Ga$-gra\-ded \ko-linear functor from
$\VectG\boticoc\VectG$ to $\VectG$ and thus endows \VectG\ with the structure of a \cgmc.

\begin{rem}
Taking natural transformations between $\Ga$-graded \ko-linear functor as 
2-mor\-phisms, \CatG\ is promoted to a 2-category. 
Moreover, the information supplied by the presence of the bicharacter $\coc$ in the
$\coc$-twisted interchange law \eqref{eq:interchangelaw4vectG} is the same the one
given by an interchanger that provides the structure of a Gray monoid on that
2-category. More specifically, the twisted interchange law is formulated for
a quadruple of morphisms and constitutes a middle step in which the composition and
the tensor product $\oticoc$ are exchanged, while the interchanger of a Gray monoid
(see e.g.\ the formula $(\mathrm D4)$ in \Cite{Rem.\,2.4}{jopR})
involves only a pair of morphisms and directly relates the two possible ways 
in which their composite can be formed. The latter description amounts
to taking two of the four morphisms in the twisted interchange law to be identitites.
None of the arguments and results in this paper relies on this 2-categorical picture
of \CatG.
\end{rem}

%%%%%%%%%%%%%%%%%%%%%%%%%%%%%%%%%%%%%%%%%%%%%%%%%%%%%%%%%%%%%%%%%%%%%%%%

\section{A tensor product of modules over a \graco\ algebra} \label{sec:oti}

Let \C\ be a braided category, $\Ga$ a finite abelian group, and $A$ a $\Ga$-graded 
algebra in \C\ that is \graco\ with respect to a two-cocycle $\coc$ on $\Ga$. Then by
the defining equality \eqref{eq:def:graco}, the opposite algebra $\Ao$ coincides 
with the cocycle-twisted algebra $A\ococ$ described in Remark \ref{rem:twist-alg}.
Now by Lemma \ref{lem:Wm}, via the braiding one obtains a right $\Ao$-module from 
a left $A$-module, while by Lemma \ref{lem:mtilde}, via cocycle-twisting one obtains
a right $A\ococ$-module from a $\Ga$-graded right $A$-module. We combine the two 
procedures by setting, for a $\Ga$-graded right $A$-module $m \eq (\dtm,\ohr)$ 
with $\dtm \eq \bigoplus_i n_i$,
  \be
  \rho\ococbr_{i,j} := \coc(i,j)\inv \, \ohr_{j,i} \cir \br_{A_i,m_j}
  \label{eq:def:invbraid+coc}
  \ee
for $i,j\iN\Ga$.

\begin{prop} \label{prop:untwisting}
Let $m \eq (\dtm,\ohr)$ be a $\Ga$-graded right module over a \graco\ algebra $A$.
 \\[2pt]
(i) We have
  \be
  m\ococbr := (\dtm,\rho\ococbr) \in A\Mod \,,
  \ee
i.e.\ $\coc\inv$-twisting the left $\Ao$-module obtained from $m$ by the
braiding yields an ordinary left $A$-module.
 \\[2pt]
(ii) 
If and only if the two-cocycle $\coc$ with respect to which $A$ is \graco\ is a 
bicharacter, the triple $(\dtm,\rho\ococbr,\ohr)$ is an $A$-bi\-module.
\end{prop}

\begin{proof}
(i)
Since, as noticed above, the opposite algebra $\Ao$ of a \graco\ algebra $A$
coincides with $A\ococ$, we can regard the left $\Ao$-module that is obtained
from $m$ via the braiding as a left $A\ococ$-module. The claim then follows by
applying Lemma \ref{lem:mtilde} `backwards' to the so obtained left $A\ococ$-module.
 \\[2pt]
(ii)
By direct calculation one finds
  \be
  \rho\ococbr_{i,j+k}  \circ (\id_{A_i} \oti \ohr_{j,k})
  = \coc(i,j)\,\coc(i,k)\,\coc\inv(i,j{+}k)\; \ohr_{i{+}j,k}
  \circ (\rho\ococbr_{i,j}\oti \id_{A_k})
  \ee
for $i,j,k\iN\Ga$. So the left and right $A$-actions on $\dtm$ commute,
$\rho\ococbr \cir (\ida \oti \ohr) \eq \ohr \cir (\rho\ococbr \oti \ida)$, if and
only if the cocycle $\coc$ is a character in its second argument. The statement then
follows because a cocycle is a character in one argument iff it is a bicharacter.
\end{proof}

\begin{rem} \label{rem:untwisting-left}
Analogously, for a $\Ga$-graded left $A$-module $m \eq (\dtm,\rho)$ with underlying
object $\dtm \eq \bigoplus_i m_i$ we set
$\ohr\ococbr_{i,j} \,{:=}\, \coc(i,j)\inv \rho_{j,i} \cir \br_{m_i,A_j}$
for $i,j\iN\Ga$. Then if $A$ is \graco, $\coc\inv$-twisting the right $\Ao$-module
obtained from $m$ by the braiding yields an ordinary right $A$-module, 
$(\dtm,\ohr\ococbr) \iN \moD A$, and iff moreover $\coc$ is a bicharacter, 
$(\dtm,\rho,\ohr\ococbr)$ is an $A$-bi\-module.
The first statement is shown in the same way as part (i) of Proposition
\ref{prop:untwisting}, while for the second statement one calculates
  \be
  \ohr\ococbr_{i+j,k} \circ (\rho_{i,j} \oti \id_{A_k})
  = \coc(i,k)\,\coc(j,k)\,\coc\inv(i{+}j,k) \, \rho_{i,j+k} \circ
  (\id_{A_i} \oti \ohr\ococbr_{j,k})
  \ee
for $i,j,k\iN\Ga$, meaning that for the left and right $A$-actions to commute,
$\coc$ now has to be a character in its first argument, which is again
equivalent to $\coc$ being a bicharacter.
\end{rem}

\begin{rem} \label{rem:graco-iff-bichar}
By (the proof of) Lemma \ref{lem:graco-iff-bichar}, the requirement that $\coc$
is a bicharacter is automatically satisfied (and in fact equivalent to $A$ being 
\graco) if the \graco\ algebra $A$ is a twisted group algebra.
\end{rem}

With the help of Proposition \ref{prop:untwisting} we obtain

\begin{prop} \label{prop:tildeoti}
Let $A$ be an algebra in a braided monoidal category \C\ for which the tensor product
$\otA$ of $A$-bimodules exists. Let $\Ga$ be a finite abelian group and let $A$ be 
\graco\ with respect to a cocycle $\coc$ on $\Ga$ that is a bicharacter on $\Ga$.
 \\[2pt]
(i) The category $\modGa A$ of $\Ga$-graded right $A$-modules and homogeneously even 
module morphisms admits a structure $(\modGa A,\totI)$ of a monoidal category.
 \\[2pt]
(ii) If $A$ is in addition \deltasep\ Frobenius, then the tensor product $m \toti n$ 
of two $\Ga$-graded right $A$-modules $(\dt m,\ohr_m)$ and $(\dtn,\ohr_n)$ 
is the image of the morphism
  \be
  p_{m,n}\ococ := (\ohr_m\oti\rho_n\ococbr)
  \circ (\id_{\dtm} \oti (\comul\cir\eta) \oti \id_{\dtn})
  \label{eq:pcocmn}
  \ee
(which is an idempotent).
\end{prop}

\begin{proof}
(i)
Consider two right $A$-modules $m \eq (\dtm,\ohr_m)$ and $n \eq (\dtn,\ohr_n)$.
Invoking Proposition \ref{prop:untwisting}, we can promote $m$ and $n$ to 
$A$-bimodules $(\dtm,\rho_m\ococbr,\ohr_m)$ and $(\dtn,\rho_n\ococbr,\ohr_n)$,
which we still denote by $m$ and $n$. The bimodule tensor product over $A$, as given
in Definition \ref{def:ota-coequ}, then yields an $A$-bimodule
  \be
  m \toti n := (\dtm,\rho_m\ococbr,\ohr_m) \ota (\dtn,\rho_n\ococbr,\ohr_n) \,.
  \ee
By forgetting the left $A$-action, this becomes
a right $A$-module, which we still denote by $m \toti n$.
Since $\otA$ is a monoidal structure on $A\bimod$, with associator and unitors inherited 
from \C, it follows immediately that $\totI$ is a monoidal structure on $\modGa A$, again 
with associator and unitors inherited from those of \C, and with $A$ as monoidal unit. 
 \\[2pt]
(ii)
If 
\C\ is abelian and
$A$ is \deltasep\ Frobenius, then the tensor product $\otA$ can be described,
as in Lemma \ref{lem:mAn=Im(p)}, as the image of an idempotent. This idempotent is
given by formula \eqref{eq:def:bimod-p}, with appropriate right and left $A$-actions.
In the case at hand the latter are $\ohr_m$ and $\rho_n\ococbr$, respectively,
thus yielding \eqref{eq:pcocmn}.
\end{proof}

\begin{rem} \label{rem:tildeoti-left}
The following analogous results hold for left modules.
For $A$ an algebra in a braided monoidal category \C\ that is \graco\ with respect
to a cocycle $\coc$ on $\Ga$, the category $A\Mod_\Ga$ of $\Ga$-graded left 
$A$-modules admits a structure $(A\Mod_\Ga,\totI)$ of a monoidal category. If 
\C\ is abelian and
$A$ is in addition \deltasep\ Frobenius, then the tensor product $m \toti n$ of
two left $A$-modules $m$ and $n$ is the image of the idempotent
$(\ohr\ococbr_m\oti\rho_n) \cir (\id_{\dtm} \oti (\comul\cir\eta) \oti \id_{\dtn})$.
In the particular case that $A$ is a commutative algebra, this reproduces results
that have been obtained, under some insignificant further conditions on \C, in
\Cite{Cor.\,1.22}{kios} and (part of) \Cite{Prop.\,3.21(i)}{ffrs}.
\end{rem}

%%%%%%%%%%%%%%%%%%%%%%%%%%%%%%%%%%%%%%%%%%%%%%%%%%%%%%%%%%%%%%%%%%%%%%%%

\section{Induced modules over twisted group algebras} \label{sec:inDA}

In Section \ref{sec:oti} we have dealt with the category $\modGa A$ of 
$\Ga$-graded right modules over a $\Ga$-gra\-ded algebra $A$; by Definition
\ref{def:Ga-module}, $\modGa A$ has homogeneously even module morphisms as morphisms. 
It is, however, equally natural to refrain from restricting to homogeneously even
morphisms. We denote by $\modeGa A$ the category having $\Ga$-graded right 
$A$-modules as objects and all module morphisms between them as morphisms.
 
In full generality, there need not be any interesting grading on the morphisms of
$\modeGa A$. However, upon induction, morphisms in \C\ give rise to homogeneously
even morphisms in $\modeGa A$, so to attain a non-trivial grading there need
to be module morphisms that are not induced. 
Henceforth we therefore restrict our attention to a situation in which it is 
guaranteed that all modules under consideration are $\Ga$-graded and that their
morphism spaces are naturally graded as well and can be described explicitly:
we require $A \eq A_{(\Ga)}$ to be a \emph{twisted group algebra} in the sense 
of Example \ref{exa:twistedgroupalg} and consider induced $A$-modules.

\begin{defi} \label{def:IndA}
An \emph{induced} (or \emph{free}) right module over an algebra $(A,\mul,\eta)$ 
in a monoidal category \C\ is a right $A$-module of the form 
$(\dtm,\ohr) \eq \II x \,{:=}\, (x \oti A,\id_x \oti \mul)$ for some $x\iN\C$.
 \\
The category $\inD A$ is the full subcategory of $\moD A$ whose objects are
induced right $A$-modules.
\end{defi}

Induced left modules and induced bimodules are defined analogously. The functor 
  \be
  \I\Colon \C \rarr~ \moD A
  \label{eq:def:Ind}
  \ee
that maps $x\iN\C$ to $\II x$ and $f \iN \Hom(x,y)$ to
$f \oti \ida \iN \HomA(\II x,\II y)$ for
$x,y\iN\C$ is called the \emph{induction} (or \emph{free}) functor.
The functor $\I$ is faithful and is left adjoint to the forgetful functor 
$U\colon \moD A \Rarr~ \C$ which acts as $(\dtm,\ohr) \Mapsto \dtm$. If $A$ is 
Frobenius (such as in the case at hand, where $A$ is a twisted group algebra), then 
$\I$ is also right adjoint to $U$, so that we have a two-sided adjunction 
$\I \,{\dashv}\, U \,{\dashv }\,\I$, i.e.
  \be
  \HomAA(\II x,m) \cong \Hom(x,\dtm) \qquand \HomAA(m,\II x) \cong \Hom(\dtm,x) \,.
  \label{eq:IUI}
  \ee

Next we make use of the fact that any $\Ga$-graded algebra $A$ is automatically 
$\Ga$-graded as a bimodule over itself. By simply assigning grade $0$ to an object 
$x \iN \C$, whereby $x \oti A$ is naturally graded as $x \oti A 
\eq \bigoplus_{i\in\Ga} x_i$ with $x_i \eq x \oti A_i$, this directly implies

\begin{lem} \label{lem:inducedGa-module}
Every induced module over a $\Ga$-graded algebra is naturally $\Ga$-graded
in the sense of Definition \ref{def:Ga-module}.
\end{lem}

As a consequence, for a twisted group algebra $A$ the category $\inD A$ is 
$\Ga$-graded in the sense of Definition \ref{def:cocgracat}.\ref{item:cocgracat}.
The $\Ga$-grading on $\inD A$ is in fact \emph{strong}, in the sense that for
every object $\II x$ and every $i\iN\Ga$ there exists an isomorphism of 
grade $i$ with domain $\II x$, namely $\II x \Rarr\cong \II{x\oti A_i}$.

By construction, all morphisms in $\inD A$ that are in the image of the
induction functor \eqref{eq:def:Ind} are homogeneously even. Non-even and 
non-homogeneous endomorphisms will be of particular interest for objects $x\iN\C$ 
for which the subset
  \be
  \stab x := \{ i \iN \Ga \,|\, x \oti A_i \Cong x \}
  \ee
of $\Ga$ is non-trivial, i.e.\ which obey if $x \oti A_i \Cong x$ for some homogeneous
component of $A$ with $i \nE 0$. We call $\stab x$ the \emph{stabilizer} of $x$.
Note that if $x \oti A_i \Cong x \Cong x \oti A_j$ for $i,j\iN\Ga$, then also
$x \oti A_{i+j} \Cong x$, i.e.\ the stabilizer is in fact a subgroup of $\Ga$.
 
\begin{exa}
Consider the semisimple modular tensor category $\C(\mathrm A_1^{(1)},k)$ 
introduced in Example \ref{exa:zz}, having non-isomorphic simple objects $U_\nu$ 
with $\nu \iN \{0,1,...\,,k\}$ and Picard group $\zz$. The stabilizers of the
simple objects with respect to the group algebra $A \eq U_0 \,{\oplus}\, U_k$ 
are trivial, $\stab{U_\nu} \eq \{0\}$, unless $\nu \eq k/2$,
in which case $\stab{U_k/2} \eq \zz$.
\end{exa}

To study the structure of morphisms in $\inD A$ in detail, we now further assume 
that \C\ is \ko-linear abelian and that the monoidal unit of \C\ is simple.
If \C\ has objects with non-trivial stabilizer, then $\inD A$ is 
not idempotent-complete and hence not abelian, 
but we can still define an object $m \iN \inD A$ to be \emph{simple} 
iff it has precisely two quotient objects, namely itself and the zero object.
We refer to an induced $A$-module that is simple as an object of $\inD A$ as
a \emph{simple induced $A$-module}, even if it would not be simple as an object
in the category of all $A$-modules.

\begin{lem}
Let \C\ be \ko-linear abelian monoidal with simple monoidal unit, let $\Ga$ be an 
abelian group, and let $A \eq A_{(\Ga)}$ be a $\Ga$-graded twisted group algebra.
An object $\II x \iN \inD A$ is a simple induced $A$-module if and only if
$x \iN \C$ is a simple object of \C.
\end{lem}

\begin{proof}
(i) If $x\iN\C$ is not simple, then it has a non-trivial quotient object, so there 
exists a non-zero $y\iN\C$, not
isomorphic to $x$, with an epimorphism $f \iN \Hom(x,y)$. We then also have
$x \oti A \,{\not\cong}\, y\oti A$ and hence $\II x \,{\not\cong}\, \II y$.
Moreover, since the functor $\I$ is a left adjoint, it preserves epimorphisms,
and hence $\II f \iN \HomA(\II x,\II y)$ is epi as well.
It follows that $\II x \iN \inD A$ is not simple.
 \\[2pt]
(ii)
Assume that the induced module $\II x$ is not simple, i.e.\ 
has a non-trivial quotient object in $\inD A$. Then there exists a non-zero
$m\iN \inD A$, not isomorphic to $\II x$, with an epimorphism in $\HomA(\II x,m)$.
By definition of $\inD A$ there is thus an object $y\iN\C$ with $m \eq \II y$.
Using that $\I$ is left adjoint to the forgetful functor $U$ it follows that
  \be
  \Hom(x,y \oti A) = \Hom(x,U\II y) \cong \HomA(\II x,m) \ne 0 \,.
  \label{Homx,yAne0}
  \ee
Assume now that $x$ is simple. If $y$ is simple as well, then it follows from
\eqref{Homx,yAne0} that $y \Cong x \oti A_i$ for some $i\iN \Ga$, and thus that
$\II x \Cong \II y \eq m$, which is a contradiction.
In case $y$ is not simple, then it has a non-trivial simple quotient $y'\iN\C$ with
an epimorphism in $\HomC(y,y')$, and thus in particular
$\HomA(\II y,\II{y'}) \Cong \Hom(y,y'\oti A) \nE 0$.
This implies, in turn, 
situation is reduced to the previous one with $y'$ in place of $y$.
\end{proof}

In a \ko-linear abelian monoidal category \C\ with simple monoidal unit
\emph{Schur's Lemma} holds (see e.g.\ \Cite{Lemma\,1.5.2}{EGno}): any non-zero 
morphism between two simple objects in \C\ is an isomorphism, and $\Hom(x,x)$ is
a \ko-division algebra for any simple object $x\iN\C$. Thus, since we assume \ko\
to be algebraically closed, for any pair of simple objects $x,y\iN\C$ we have
$\Hom(x,y) \Cong \ko$ if $x$ and $y$ are isomorphic, and $\Hom(x,y) \eq 0$ else.
We will show that in this situation there is a graded analogue of Schur's Lemma
for induced modules over twisted group algebras in \C. We first observe

\begin{prop} \label{prop:dimEndAA=ko.stab}
Let $A$ be a $\Ga$-graded twisted group algebra in a \ko-linear abelian 
monoidal category \C\ with simple monoidal unit, with $\Ga$ a finite abelian group.
Let $x\iN \C$ be simple, with stabilizer $\stab x$. Then there is an isomorphism
  \be
  \EndAA(\II x) \cong \kotw \stab x 
  \label{EndAA=k.stab}
  \ee
of $\stab x$-graded \ko-algebras between the endomorphims of the induced module 
$\II x \iN \inD A$ and the group algebra $\ko\, \stab x$ twisted by a
two-cocycle $\cocs$.
\end{prop}

\begin{proof} 
We first establish \eqref{EndAA=k.stab} as an isomorphism of vector spaces.
Abbreviate $\stab x \,{\equiv}\, \sta$.
Select for each $\sta$-coset $c$ of $\Ga$ a representative $i_c \iN \Ga$, so that
$c \eq i_c\pl\sta \eq \{ i_c \pl s \,|\, s\iN \sta \}$. Then as a set we
have $\Ga \eq \{ i_c \pl s \,|\, c \iN \Ga/\sta,\, s\iN \sta \}$ and thus
  \be
  x\oti A =\, \bigoplus_{s \in \sta} \bigoplus_{c\in \Ga/\sta}\! x \oti A_{i_c+s}
  \cong\, \bigoplus_{s \in \sta} \Big(\! \bigoplus_{c\in \Ga/\sta}\!
  x \oti A_s \oti A_{i_c} \Big)
  \cong\, \bigoplus_{s \in \sta} \Big(\! \bigoplus_{c\in \Ga/\sta}\!
  x \oti A_{i_c} \Big) .
  \label{eq:xA-decompoH}
  \ee
In short, $x\oti A$ decomposes into a direct sum of $|\sta|$ copies of the 
`$\Ga$-orbit' of $x$ which is the object $\bigoplus_{c\in \Ga/\sta} x \oti A_{i_c}$,
where up to isomorphism the choice of representatives $i_c$
of the cosets $c \iN \Ga/\sta$ is irrelevant.
It follows that as vector spaces we have
  \be
  \begin{aligned}
  \End_{\!A}(\II x) &
  \eqc{eq:IUI} \Hom(U(\II x),x) = \HomC(x\oti A,x)
  \\
  & \eqc{eq:xA-decompoH} \bigoplus_{s \in \sta} \bigoplus_{c\in \Ga/\sta}\!\!
  \Hom(x\oti A_{i_c},x)
  \cong\, \bigoplus_{s \in \sta} \Hom(x,x) \,\cong\, \ko^{|\sta|} .
  \end{aligned}
  \label{eq:EndAIndx=kH}
  \ee
To see that $\End_{\!A}(\II x)$ is isomorphic to $\ko\sta$ as a $\sta$-graded
vector space, we provide a graded basis. For $s\iN\sta$ select an isomorphism
$\varphi_s \iN \Hom(x,x\oti A_s)$ and define an endomorphism
$\phi_s \,{:=}\, \bigoplus_{i\in\Ga} \phi_{s;i}$ of $x\oti A$ by
  \be
  \phi_{s;i}\Colon x\oti A_i \rarr{\varphi_s\oti\id}
  (x\oti A_s) \oti A_i \rarr{\ass_{x,A_s,A_i}} x\oti (A_s \oti A_i)
  \rarr{\id \oti \mul_{s,i}} x\oti A_{i+s} \,.
  \ee
It follows directly from associativity of $A$ that the morphisms $\phi_s$ are in 
fact module endomorphisms, $\phi_s \iN \EndAA(\II x)$. Also, $\phi_s$ is homogeneous
of grade $s$, so that in particular the set $\{ \phi_s \,|\, s\iN\sta \}$ is linearly
independent and is thus a graded basis of $\EndAA(\II x) \Cong \ko^{|\sta|}$.
Finally we have $\phi_{s'}\cir\phi_s \eq \cocs(s,s')\,\phi_{s+s'}$ with
numbers $\cocs(s,s')\iN\kox$ for $s,s'\iN\sta$. 
Associativity of composition amounts to the statement that the numbers
$\cocs(s,s')$ assemble to a two-cocycle on $\sta$.
Modifying the choice of isomorphisms $\varphi_s$ changes $\cocs$ by a two-coboundary.
\end{proof} 

\begin{thm}(\emph{\GSL}) \label{thm:GSL}
 \\
Let $A$ be a $\Ga$-graded twisted group algebra in a \ko-li\-ne\-ar abelian  
monoidal category \C\ with simple monoidal unit, with $\Ga$ a finite abelian group.
\\[2pt]
(i) The space of module morphisms between any two simple induced $A$-modules $\II x$
and $\II y$ is either zero or isomorphic as a $\Ga$-graded \ko-vector space to the 
spaces of endomorphisms of both $\II x$ and $\II y$.
 \\[2pt]
(ii) Any non-zero homogeneous module morphism between simple induced $A$-modules
is an isomorphism.
\end{thm}

\begin{proof} 
(i) Analogously as in \eqref{eq:EndAIndx=kH} we have
$\HomA(\II x,\II y) \Cong \HomC(x\oti A,y)$. 
If $x \,{\not\cong}\, y\oti A_i$ for every $i\iN\Ga$, then this space is zero.
Otherwise we have $x \Cong y\oti A_j$ for some $j\iN\Ga$ (which implies that
$\stab y \eq \stab x$ and thus $\End_{\!A}(\II y) \Cong \End_{\!A}(\II x)\:$), and
again analogously as in \eqref{eq:EndAIndx=kH} we obtain
  \be
  \HomA(\II x,\II y) \cong \bigoplus_{s \in \stab x} \bigoplus_{c\in \Ga/\stab x}\!\!
  \Hom(x\oti A_{i_c+j},x) \,.
  \label{eq:HomAIIxIIy}
  \ee
Now for any $j\iN\Ga$ there is a permutation $\pi$ of $\Ga/\stab x$ such that
for every $c\iN\Ga/\stab x$ one has $i_c\pl j \eq i_{\pi(c)}\pl s'$ with some
$s'\iN\stab x$. This implies that the right hand side of \eqref{eq:HomAIIxIIy}
is isomorphic to 
$\bigoplus_{s \in \stab x}\! \bigoplus_{c\in \Ga/\stab x}\! \Hom(x\oti A_{i_c},x)$,
and hence $\HomA(\II x,\II y) \Cong \End_{\!A}(\II x)$ as vector spaces. It 
then also follows directly
that $\HomA(\II x,\II y)$ inherits a $\Ga$-grading from $\EndAA(\II x)$.
 \\[2pt]
(ii) Via composition, the space $\HomA(\II x,\II y)$ has a natural structure of an
$\EndAA(\II x)\text-
            $\linebreak[0]$
\EndAA(\II y)$-bimodule. As a consequence the homogeneous subspaces of the 
$\Ga$-graded vector space $\HomA(\II x,\II y)$ are one-dimensional.
\end{proof} 

\begin{rem}
As mentioned, for abelian monoidal categories the ordinary Schur Lemma is valid.
Accordingly, if the category $\inD A$ is abelian, then Theorem \ref{thm:GSL} implies
the ordinary Schur Lemma. But generically, $\inD A$ is not idempotent complete and 
thus not abelian. On the other hand, at least if $A$ is a special Frobenius algebra,
every $A$-module is a module retract of an induced module \Cite{Lemma\,4.8}{ffrs}
and the category $\moD A$ of all $A$-modules is the idempotent completion of $\inD A$.
As a consequence, if \C\ is abelian and $A$ is special Frobenius, then so is $\moD A$,
so that the ordinary Schur Lemma holds in $\moD A$.
Also, as follows directly from the proof of Theorem \ref{thm:GSL}, the ordinary 
Schur Lemma holds in the wide subcategory $\inDe A$ of $\inD A$ whose morphisms
are the homogeneously even morphisms of $\inD A$.
\end{rem}

\begin{exa}
For $\Ga \eq \zz$, i.e.\ $A \eq \one \,{\oplus}\, L$ with $L$ an order-2
invertible object, we have
  \be
  \HomA(\II x,\II y) \,\cong\, \left\{\!\! \begin{array}{cl}
  \ko & \text{if}~ x\oti L \,{\not\cong}\, x ~\text{and}~
        (x \Cong y ~\text{or}~ x \Cong y \oti L\,) \,,
  \nxl2
  \ko(\zz) & \text{if}~ x\oti L \Cong x \Cong y \,,
  \nxl2 
  0 & \text{otherwise} \,. \eear\right.
  \label{eq:210}
  \ee
\end{exa}

\begin{rem} \label{rem:SSL}
Particular realizations of Proposition \ref{prop:dimEndAA=ko.stab} and Theorem
\ref{thm:GSL}, or of parts thereof, are well known, like e.g.\ Theorem 14 of
\cite{maZh3} which describes homogeneous automorphisms of graded modules over
graded Lie algebras. In the representation theory of superalgebras,
the trichotomy \eqref{eq:210}, then sometimes called the \emph{Super Schur Lemma}, 
lists the possible dimensions of morphism spaces between simple modules over
a finite-dimensional associative superalgebra, see e.g.\
Corollary 2.15 of \cite{joze} and Lemma 2.3 of \cite{waWan}.
In that case, the possibility of having stabilizer $\zz$ is coupled to the
existence of the simple superalgebras of type $Q(n)$; accordingly, modules with
stabilizer $\zz$ are then called modules of type Q.
In the context of spin topological quantum field theories and their applications 
to fermionic topological phases in condensed matter physics, objects with 
stabilizer $\zz$ are also known as \emph{Majorana objects}, see e.g.\ \cite{aalW}.
It is, however, worth keeping in mind that what is relevant for the possible
existence of such objects is only that the relevant algebra $A$ is \zz-graded, not 
that it admits any kind of super structure (compare also Example \ref{exa:zz-2}).
\end{rem}

%%%%%%%%%%%%%%%%%%%%%%%%%%%%%%%%%%%%%%%%%%%%%%%%%%%%%%%%%%%%%%%%%%%%%%%%

\section{A graded-monoidal structure for induced modules} \label{sec:gil}

In this section we take \C\ to be braided and $A$ to be an algebra in \C\ that is
\graco\ for some finite abelian group $\Ga$ and normalized two-cocycle $\coc$. 

As seen above, the category $\inD A$ is naturally $\Ga$-graded. Also, as a full
subcategory of $\modGa A$, the category $\inDe A$ inherits the tensor product $\totI$
that we considered in Section \ref{sec:oti} (see Proposition \ref{prop:tildeoti}).
Indeed, the braiding of \C\ gives a (non-canonical) isomorphism $U\II x \oti U\II y 
\eq x \oti A \oti y \oti A \rarr\cong x \oti y \oti A \oti A \eq U\II{x\oti y} \oti A$.
This isomorphism fits together with the $A$-action on $\II x$ and $\II y$ in such a 
way that the coequalizer in the definition of $\II x \ota \II y$ becomes isomorphic 
as an object to $x \oti y \oti B$ with $B$ the coequalizer of $\ida \oti \mul$ and
$\mul \oti \ida$ on $A\oti A\oti A$, i.e.\ $B \Cong A\ota A \Cong A$; moreover, 
the so obtained isomorphism from $U(\II x \ota \II y)$
to $x \oti y \oti A \eq U\II{x\oti y}$ is compatible with the $A$-action
on $\II{x\oti y}$.
	
An obvious question is then whether the monoidal structure $\totI$ on $\inDe A$ 
can be extended to all 
of $\inD A$. As we will see, this is the case only if $A$ is in fact commutative.
In the general case it is natural to try to construct instead a suitably twisted
monoidal structure on $\inD A$ that makes it a \cgmc\ in
the sense of Definition \ref{def:cocgramoncat}, and that reduces to the tensor 
product $\totI$ for homogeneous morphisms of grade zero. 

To find a candidate for such a twisted monoidal structure on $\inD A$, it will be
convenient to consider the monad $\AA$ with underlying endofunctor $T_A\eq {-}\oti A$,
as described in \eqref{eq:def:AA}, and to identify a corresponding candidate for a
twisted monoidal structure on $\AA$. Working with $\AA$, we realize $\inD A$ as 
the Kleisli category \Cite{Ch.\,4.1}{BOrc2} over $\AA$.

\begin{defi} \label{def:Kleisli}
The \emph{Kleisli category} over a monad $\TT \eq (T,\mulT,\etaT)$ on a category
\C\ is the following category $\C_\TT$: The objects of $\C_\TT$ are the same as 
those of \C\ -- for distinction, we write $\ol x$ for the object $x\iN\C$ when 
regarded as an object of $\C_\TT$. The morphisms of $\C_\TT$ are given by
$\Hom_{\C_\TT}(\ol x,\ol y) \eq \HomC(x,Ty)$.
The composition $\Hom_{\C_\TT}(\ol y,\ol z) \Times \Hom_{\C_\TT}(\ol x,\ol y)
\Rarr~ \Hom_{\C_\TT}(\ol x,\ol z)$ of morphisms is the \emph{Kleisli composition}
  \be
  (g,f) \,\xmapsto{~~~}\, \mulT_z \circ Tg \circ f =: g \ocir f
  \label{eq:Kleislicomp}
  \ee
for $f \iN \Hom_{\C_\TT}(\ol x,\ol y)$ and $g \iN \Hom_{\C_\TT}(\ol y,\ol z)$.
The identity morphism of $\ol x$ is the relevant component of the unit of the monad,
$\id^{}_{\ol x} \eq \etaT_x \colon x \Rarr~ Tx$.
\end{defi}

Analogously as for algebras, an \emph{induced} (or free) module over a monad $\TT$
on \C\ is a $\TT$-mo\-du\-le of the form $\IIT x \,{:=}\, (Tx,\mulT_x)$ for some 
$x \iN \C$.
Of significance to us is that the category $\TT\Ind$ of induced modules over a monad 
$\TT$ is equivalent to the Kleisli category $\C_\TT$ of $\TT$. An equivalence
$\C_\TT \Rarr\simeq \TT\Ind$ is given by mapping objects as 
$\ol x \,{\xmapsto{~~}}\, \IIT x$ and morphisms as 
$(\ol x \Rarr f \ol y) \,{\xmapsto{~~}}\, \mulT_y \cir Tf$ \Cite{Prop.\,4.1.6}{BOrc2};
as a quasi-inverse this has the functor given by $\IIT x \,{\xmapsto{~~}}\, \ol x$ 
and $(\II x \Rarr g \II y) \,{\xmapsto{~~}}\, g \cir \etaT_x$.
  
Furthermore, in the case of the monad $\AA$, the category $\inD A$ of induced
$A$-modules is isomorphic to the category of induced $\AA$-modules 
and thus equivalent to the Kleisli category $\C_\AA$.
The Kleisli composition \eqref{eq:Kleislicomp} on $\C_\AA$ reads
  \be
  g \ocir f = (\id_z\oti\mul) \circ (g\oti\ida) \circ f \,.
  \ee

A monoidal structure $\tZ$ on a monad $\TT$ on a monoidal category \C\ 
equips its Kleisli category $\C_\TT$ with a tensor product $\otK$ via setting 
$\ol x \otk \ol y \,{:=}\, \ol{x\oti y}$ for $x,y\iN\C$ and
  \be
  f \otk g := \TZ{x'}{y'} \cir (f \oti g) \Colon x\oti y \rarr~ T(x'{\otimes}y')
  \label{eq:def:otk}
  \ee
for morphisms $f\colon x\Rarr~Tx'$ and $g\colon y\Rarr~Ty'$ in \C.
Indeed, there is a bijection between monoidal structures on $\TT$ and monoidal 
structures on $\C_\TT$ for which the induction functor
from \C\ to $\C_\TT$ is strict monoidal \Cite{Prop.\,1.2.2}{seal}.
Thus if $\TT$ is endowed with a monoidal structure, then there is a monoidal 
structure on the category of induced $\TT$-modules such that the equivalence
$\C_\TT \Rarr\simeq \TT\Ind$ is strict monoidal.

Now recall from Definition \ref{def:monoidalmonad} that
a monoidal structure on a monad requires the multiplication of the monad
to become a monoidal natural transformation. We observe

\begin{lem} \label{lem:Kleisli-exclaw}
Let $\tZ$ be a monoidal structure on a monad $\TT \eq (T,\mulT,\etaT)$. The
monoidality of the multiplication $\mulT$ is equivalent to the validity of the
interchange law
  \be
  (f' \otk g') \ocir (f \otk g) = (f' \ocir f) \otk (g' \ocir g)
  \label{eq:Kleisliinterchange}
  \ee
for any quadruple of morphisms $f\colon x\Rarr~ Tx'$, $g\colon y\Rarr~ Ty'$,
$f'\colon x'\Rarr~ Tx''$ and $g'\colon y'\Rarr~ Ty''$ in \C.
\end{lem}

\begin{proof}
Monoidality of $\mulT$ is expressed by commutativity of the diagram 
\eqref{eq:mulT.montrafo}, i.e.\ by the equality
$\mulT_{x\otimes y} \cir T(\TZ xy) \cir \TZ{Tx}{Ty}
\eq \TZ xy \cir (\mulT_x \oti \mulT_y)$. Precomposing with 
$(Tf' \oti Tg') \cir (f \oti g)$, the left hand side of this equality becomes,
using naturality of $\tZ$,
  \be
  \begin{aligned}
  \mulT_{x''\otimes y''} \circ T(\TZ{x''}{y''})
  & \circ \TZ{Tx''}{Ty''} \!\circ (Tf' \oti Tg') \circ (f \oti g)
  \nxl1
  & = \mulT_{x''\otimes y''} \cir T\big( \TZ{x''}{y''} {\circ}\, (f' \oti g') \big)
  \circ \big( \TZ{x'}{y'} {\circ}\, (f \oti g) \big)
  \equiv (f' \otk g') \ocir (f \otk g) \,,
  \end{aligned}
  \label{eq:Kleisli-exclaw-lhs}
  \ee
while the right hand side becomes
  \be
  \begin{aligned}
  \TZ{x''}{y''} \circ (\mulT_{x''} \oti \mulT_{y''})
  & \circ (Tf' \oti Tg') \circ (f \oti g)
  \nxl1
  & = \TZ{x''}{y''} \circ \big( (\mulT_{x''} \cir Tf' \cir f)
  \oti (\mulT_{y''} \cir Tg' \cir g) \big)
  \equiv (f' \ocir f) \otk (g' \ocir g) \,.
  \end{aligned}
  \label{eq:Kleisli-exclaw-rhs}
  \ee
Thus commutativity of \eqref{eq:mulT.montrafo} implies \eqref{eq:Kleisliinterchange}.
Conversely, when taking $x \eq Tx' \eq T^2x''$ and $y \eq Ty' \eq T^2y''$ as well as
$f \eq \id_{Tx'}$ etc.,
\eqref{eq:Kleisliinterchange} reduces to commutativity of \eqref{eq:mulT.montrafo}.
\end{proof}

\begin{rem}
In case $A$ is a \deltasep\ Frobenius algebra (e.g.\ a twisted group algebra)
in an abelian category,
the tensor product of $A$-bimodules, and hence also $\totI$, can be expressed as the
image of the idempotent \eqref{eq:pcocmn}. That $\II x \toti \II y$ can be identified
with $\II{x{\otimes}y}$ is then understood directly by noticing that the object
underlying the module $\II x \toti \II y$ is the image of 
  \be
  \id_x \otimes \big[
  (\br\inv_{\aA,y} \oti \ida) \cir (\id_y \oti p_{A,A}^{}) \cir (\br^{}_{\aA,y} \oti \ida)
  \big] ~\iN \End(x \oti A \oti y \oti A) \,,
  \ee
with $p_{A,A}^{}$ the idempotent
  \be
  p_{A,A}^{} := (\mul \oti \mul) \circ \big(\ida \oti (\comul\cir\eta) \oti \ida \big)
  = \comul \circ \mul
  \label{eq:pAA}
  \ee
in $\End(A \oti A)$, by which $A$ is canonically exhibited as a retract of $A\oti A$. 
Also, if $A$ is separable Frobenius, then by Lemma 
\ref{lem:A-AA-Frobenius} the monad $\AA$ is separable Frobenius in the sense of 
Definition \ref{def:frobmonad}, with comonad structure given by \eqref{eq:def:coAA}.
\end{rem}

We take the explicit form of the idempotent \eqref{eq:pAA} a motivation to make the
ansatz
  \be
  \TZ xy := (\id_x \oti \id_y \oti \mul) \circ (\id_x \oti \br_{A,y} \oti \ida)
  \label{eq:def:TZ4AA}
  \ee
for $x,y\iN\C$ as a candidate for a twisted analogue of a monoidal structure on 
the monad $\AA$, even for algebras $A$ that are not Frobenius. This ansatz is 
vindicated by

\begin{lem} \label{lem:monmonad4commFrob}
 ~\\[2pt]
(i)
For $A \eq (A,\mul,\eta)$ an algebra in a braided monoidal category
$(\C,\otimes,\one;\br)$, the family \eqref{eq:def:TZ4AA} of morphisms furnishes a
monoidal structure on the monad $\AA$ if and only if $A$ is commutative.
 \\[2pt]
(ii)
If $A \eq (A,\mul,\eta,\comul,\eps)$ is a commutative \deltasep\ Frobenius algebra,
then the monoidal structure \eqref{eq:def:TZ4AA} on
$\AA$ amounts to a tensor product on the Kleisli category $\C_\AA$ that is 
isomorphic to the tensor product $\totI$ of induced $A$-modules.
\end{lem}

\begin{proof}
(i)
That the family \eqref{eq:def:TZ4AA} of morphisms gives a lax monoidal structure 
on the functor ${-} \oti A$
follows immediately from the definitions, while naturality of the family 
\eqref{eq:def:TZ4AA} and monoidality of the natural transformation $\etaA$ are 
direct consequences of the naturality of the braiding $\br$. 
Finally, that the natural transformation $\mulA\colon \AA^2\Rarr~\AA$ satisfies
\eqref{eq:mulT.montrafo} boils down to the equality
$\mul \cir (\mul \oti \mul) \cir (\ida \oti \br_{A,A} \oti \ida)
\eq \mul \cir (\mul \oti \mul)$. This is satisfied iff $A$ is commutative. 
 \\[2pt]
(ii)
By \deltasepy\ of $A$, the family $\SZ xy \,{:=}\, 
(\id_x \oti \br_{\!A,y}\inv \oti \ida) \cir (\id_x \oti \id_y \oti \comul)$ of
morphisms is a one-sided inverse to the family \eqref{eq:def:TZ4AA}, 
$\tZ \cir \sZ \eq \id$. Accordingly we obtain a tensor product on the category
$\C_\AA$ corresponding to \eqref{eq:def:TZ4AA} that
can be described as the image of the idempotent $\SZ xy \cir \TZ xy$.
Moreover, the Frobenius property and commutativity of $A$ imply that
$\SZ xy \cir \TZ xy \eq p_{\II x,\II y}$ for all $x,y\iN\C$, with $p_{m,n}$ as
defined in \eqref{eq:def:bimod-p}, i.e.\ $\SZ xy \cir \TZ xy$ equals the idempotent
whose image is the tensor product $\II x \ota \II y$ of induced $A$-modules.
\end{proof}

According to Lemma \ref {lem:monmonad4commFrob} the family \eqref{eq:def:TZ4AA} is 
no longer a monoidal structure on $\AA$ when the algebra $A$ is \graco\ for a 
non-trivial two-cocycle $\coc$. Instead we have

\begin{prop} \label{prop:twistedinterchange}
Let $A$ be a \graco\ algebra in a \ko-linear braided monoidal category \C. Define
$\otk$ as in \eqref{eq:def:otk}, with $\tZ$ given by \eqref{eq:def:TZ4AA}. Then in 
the Kleisli category $\C_\AA$ the \emph{$\coc$-twisted interchange law}
  \be
  (f' \ocir f) \otk (g' \ocir g) 
  = \sum_{i,j\in\Ga} \coc(i,j)\, (f'_{} \otk g'_j) \ocir (f_i^{} \otk g)
  \label{eq:graco-interchange} 
  \ee
holds for all quadruples $f\colon x\Rarr~ x'\oti A$, $g\colon y\Rarr~ y'\oti A$,
$f'\colon x'\Rarr~ x''\oti A$ and $g'\colon y'\Rarr~ y''\oti A$ of morphisms, where
$f \eq \sum_{i\in\Ga\!}f_i$ and $g'\eq\sum_{j\in\Ga\!}g'_j$ are the decompositions of
$f$ and $g'$ into sums of homogeneous morphisms $f_i\colon x\Rarr~ x'\oti A_i$
and  $g'_j\colon y'\Rarr~ y''\oti A_j$, respectively.
\end{prop}

\begin{proof}
Instead of writing out the relevant formulas, which are a bit lengthy, we work with 
string diagrams. Using associativity of $A$ we have
 \def\loch{-0.5}  %  bottom of straight line
 \def\locH{3.5}   %  height straight line
 \def\locHH{4.0}  %  height longer straight line == \locH + \loch
 \def\locX{0}     %  Position of x line
 \def\locY{1.5}   %  Position of y line
  \be
  (f' \otK g') \ocir (f \otK g) ~=~~
  \raisebox{-5.6em}{ \begin{tikzpicture}[scale=\locscale]
  \drawMod (\locX,\loch) node[below=-1pt] {$x$} -- +(0,\locHH) node[above=-1pt] {$x''$};
  \drawMod (\locY,\loch) node[below=-1pt] {$y$} -- +(0,\locHH) node[above=-1pt] {$y''$};
  \erase   (\locX,1.95) to[out=90,in=200] (\locY+0.25,2.85) ;
  \drawAlg (\locX,1.95) \boxlabel{f'} to[out=90,in=200] (\locY+0.25,2.85) ;
  \erase (\locY-0.3,1.17) to[out=10,in=200] (\locY+0.55,1.35);
  \drawAlg (\locX,0.3) \boxlabel{f} to[out=90,in=200] (\locY+0.55,1.35);
  \drawAlg (\locY,1.95) \boxlabel{g'} \upto +(0.25,0.85)
           \rightmult{0.25} -- (\locY+0.5,\locH) node[above=-1pt] {$A$};
  \drawAlg (\locY,0.3) \boxlabel{g^{\phantom:}\!} \upto +(0.55,0.95) \upto (\locY + 0.75,2.85);
  \end{tikzpicture}
  }
  ~= \quad
  \raisebox{-5.6em}{ \begin{tikzpicture}[scale=\locscale]
  \drawMod (\locX,\loch) node[below=-1pt] {$x$} -- +(0,\locHH) node[above=-1pt] {$x''$};
  \drawMod (\locY,\loch) node[below=-1pt] {$y$} -- +(0,\locHH) node[above=-1pt] {$y''$};
  \erase (\locY+0.2,1.8)\upto (\locY+0.5,2.5);
  \drawAlg (\locY,1.4) \boxlabel{g'} \upto (\locY+0.5,2.5) \leftturn{0.25} \rightmult{0.25} -- (\locY+0.5,\locH) node[above=-1pt] {$A$};
  \drawAlg (\locY,0.3) \boxlabel{g^{\phantom:}\!} \upto +(0.75,1) \upto  (\locY+0.75,2.75);
  \erase (\locY-0.3,2.65) to[out=10,in=190] (\locY+0.25,2.75) ;
  \drawAlg (\locX,1.75) \boxlabel{f'} to[out=90,in=190] (\locY+0.25,2.75) ;
  \erase (0.5*\locX+0.5*\locY,1.65) to[out=40,in=190] (\locY+0.75,2.3);
  \drawAlg (\locX,0.3) \boxlabel{f}  to[out=90,in=190] (\locY+0.75,2.3);
  \end{tikzpicture}
  }
  \quad = \quad
  \raisebox{-5.6em}{ \begin{tikzpicture}[scale=\locscale]
  \drawMod (\locX,\loch) node[below=-1pt] {$x$} -- +(0,\locHH) node[above=-1pt] {$x''$};
  \drawMod (\locY,\loch) node[below=-1pt] {$y$} -- +(0,\locHH) node[above=-1pt] {$y''$};
  \drawAlg (\locY,1.3) \boxlabel{g'} \upto (\locY+0.25,2.25)
           \rightmult{0.25}\rightmult{0.25} \leftturn{0.25};
  \drawAlg (\locY,0.25) \boxlabel{g^{\phantom:}\!} \upto +(0.85,1) \upto  (\locY+1,2.5);
  \erase (\locY-0.3,2.88) to[out=30,in=190] (\locY+0.5,3);
  \drawAlg (\locX,1.75) \boxlabel{f'} to[out=90,in=180] (\locY+0.5,3) -- (\locY+0.5,\locH) node[above=-1pt] {$A$};
  \erase (0.5*\locX+0.5*\locY,1.8) to[out=50,in=190] (1.5,2) to[out=10,in=270] (\locY+0.75,2.25);
  \drawAlg (\locX,0.3) \boxlabel{f} to[out=90,in=189] (1.5,2) to[out=10,in=270] (\locY+0.75,2.25);
  \end{tikzpicture}
  }
  \label{eq:otk-ocir-otk}
  \ee
and
 \\[-2.1em]
  \be
  (f'\ocir f) \otk (g'\ocir g)
  ~=\quad
  \raisebox{-5.6em}{ \begin{tikzpicture}[scale=\locscale]
  \drawMod (\locX,\loch) node[below=-1pt] {$x$} -- +(0,\locHH) node[above=-1pt] {$x''$};
  \drawMod (\locY,\loch) node[below=-1pt] {$y$} -- +(0,\locHH) node[above=-1pt] {$y''$};
  \erase (\locY-0.3,2.85) to[out=10,in=180] (\locY+0.5,2.95);
  \drawAlg (\locX,1.5) \boxlabel{f'} \rightmult{0.3} to[out=90,in=180] (\locY+0.5,2.95);
  \drawAlg (\locX,0.3) \boxlabel{f} \upto (\locX+0.8,1.82);
  \drawAlg (\locY,1.5) \boxlabel{g'} \rightmult{0.3}
  \upto (\locY+0.75,2.7) \leftturn{0.25} \upto (\locY+0.5,\locH) node[above=-1pt] {$A$};
  \drawAlg (\locY,0.3) \boxlabel{g^{\phantom:}\!} \upto (\locY+0.8,1.82);
  \end{tikzpicture}
  }
  ~=\quad
  \raisebox{-5.6em}{ \begin{tikzpicture}[scale=\locscale]
  \drawMod (\locX,\loch) node[below=-1pt] {$x$} -- +(0,\locHH) node[above=-1pt] {$x''$};
  \drawMod (\locY,\loch) node[below=-1pt] {$y$} -- +(0,\locHH) node[above=-1pt] {$y''$}; 
  \erase (\locY-0.3,2.97) to[out=6,in=180] (\locY+0.3,2.97);
  \drawAlg (\locX,1.75) \boxlabel{f'} to[out=90,in=180] (\locY+0.5,3) -- (\locY+0.5,\locH) node[above=-1pt] {$A$};
  \erase (\locY-0.3,1.96) to[out=9,in=180] (\locY+0.5,2.02);
  \drawAlg (\locX,0.3) \boxlabel{f} to[out=90,in=180] (1.5,2) to[out=0,in=270] (\locY+0.25,2.25);
  \drawAlg (\locY,1.4) \boxlabel{g'} \upto (\locY+0.75,2.25)
           \leftmult{0.25} \rightmult{0.25} \leftturn{0.25};
  \drawAlg (\locY,0.3) \boxlabel{g^{\phantom:}\!} \upto +(0.75,1) \upto  (\locY+1,2.5);
  \end{tikzpicture}
  }
  \label{eq:ocir-otk-ocir}
  \ee
If $A$ is \graco, i.e.\ satisfies \eqref{eq:def:graco}, then the right hand sides of
\ref{eq:otk-ocir-otk} and \ref{eq:ocir-otk-ocir},
when restricted to the components $f_i$ and $g'_j$, differ by a factor $\coc(i,j)$,
thus proving the equality \eqref{eq:graco-interchange}.
\end{proof}

\begin{rem}
If $\coc$ is the trivial cocycle -- and thus in particular if $\Ga$ is trivial,
i.e.\ $A$ is commutative -- we obtain the ordinary interchange law and hence
an ordinary monoidal structure on the category $\inD A$.
\end{rem}

Now recall the notion of a \cgmc\ from Definition
\ref{def:cocgramoncat}. The previous observations can be combined to

\begin{thm} \label{thm:II}
Let $\Ga$ be a finite abelian group and $\coc$ a $\Ga$-bicharacter. Let \C\ be a 
\ko-linear brai\-ded monoidal category and $A$ a $\Ga$-graded
algebra in \C\ that is \graco\ with respect to $\coc$.
 \\[2pt]
(i)
The category $\inD A$ can be endowed with the structure of a \cgmc, with $\coc$-\twi\
monoidal structure obtained by transporting $\otk$ from $\CAi$ to $\inD A$.
 \\[2pt]
(ii)
If in addition the tensor product $\otA$ of $A$-bimodules exists and $A$ is 
\deltasep\ Frobenius, then the restriction of $\otk$ to $\inDe A$ is isomorphic to 
the tensor product $\totI$ on $\inDe A$ that is induced from the bimodule 
tensor product $\otA$.
\end{thm}

\begin{proof}
The $\coc$-twisted interchange law \eqref{eq:graco-interchange} amounts to 
a structure of a \cgmc\ $(\CA,\otk)$ on the Kleisli category of \AA. This can be
transported to $\AA\Ind$ and thus to $\inD A$ in the same way as a monoidal structure.
Claim (ii) follows directly from Lemma \ref{lem:monmonad4commFrob}(ii).
\end{proof}

\begin{rem} \label{rem:A0=1}
If the neutral component $A_0$ of $A$ is isomorphic to the tensor unit $\one$ (which
is e.g.\ the case for twisted group algebras), then every homogeneously even morphism
in $\inD A$ is in the image of the induction functor $\I\colon \C\Rarr~ \inD A$. As
a consequence, $\I$ restricts to an equivalence from \C\ to the wide subcategory 
$\inDe A$ of $\inD A$. Moreover, the restriction of the $\coc$-twisted tensor 
product on $\inD A$ to $\inDe A$ defines a tensor product on $\inDe A$ and, 
provided that the tensor product of $A$-bimodules exists,
this coincides with the tensor 
product $\totI$ that $\inDe A$ inherits as a full subcategory of $\modGa A$.
Accordingly we could equivalently have \emph{defined} the tensor product on $\inDe A$
by setting $\II x \,{\otimes_{\inDe A}}\,\II y \,{:=}\, \II{x\,{\otimes_\C}\,y}$.
\end{rem}

\begin{rem}
One notable inspiration for this text is the study of supercategories in \cite{brEl}.
The 1-categorical structure of a monoidal supercategory in \cite[Def.\,1.4(i)]{brEl}
coincides with the tensor product with which our construction would endow the
$\coc$-graded right modules of a superalgebra, with $\coc(i,j) \eq (-1)^{ij}$
the non-trivial \zz-bicharacter. That we work with right rather than left modules 
has its origin in the fact that we then produce the same super interchange law
$(f'\oticoc g') \cir (f\oticoc g) \eq (-1)^{|g'|\,|f|} (f'\cir f) \oticoc (g'\cir g)$
as in \cite{brEl}. Had we used left modules instead, we would have arrived at
$(f'\oticoc g') \cir (f\oticoc g) \eq (-1)^{|f'|\,|g|} (f'\cir f) \oticoc (g'\cir g)$,
as the braiding involved would affect the complementary components.
 \\
In the event that $A \eq A_0 \,{\oplus}\, A_1$ is a twisted group algebra over \zz,
we also get what in \cite{brEl} is called a $\Pi$-structure: it is given by tensoring
with the invertible object $\II{A_1}$. The latter is an \emph{odd unit}, i.e.\ 
odd-isomorphic to the monoidal unit $\II{A_0}$. Thus in this case $\inD A$ is a
\emph{monoidal $\Pi$-supercategory} in the sense of Definition 1.12 on \cite{brEl}.
\end{rem}

\vskip 2em

{\sc Acknowledgments:}\\[.3em]
We thank Lukas Woike for helpful comments.
J.F.\ is supported by VR under project no.\ 2022-02931.

%%%%%%%%%%%%%%%%%%%%%%%%%%%%%%%%%%%%%%%%%%%%%%%%%%%%%%%%%%%%%%%%%%%%%%%%

\vskip 3em

 \newcommand\wb{\,\linebreak[0]} \def\wB {$\,$\wb}

 \newcommand\Bi[2]    {\bibitem[#2]{#1}}
 \newcommand\J[7]     {{\em #7}, {#1} {#2} ({#3}) {#4--#5} {{\tt [#6]}}}
 \newcommand\JO[6]    {{\em #6}, {#1} {#2} ({#3}) {#4--#5} }
 \newcommand\BOOK[4]  {{\em #1\/} ({#2}, {#3} {#4})}

 \newcommand\Prep[2]  {{\em #2}, pre\-print {\tt #1}}

   \def\adma  {Adv.\wb Math.}
   \def\apal  {Ann.\wB Pure\wB Appl.\wB Logic}
   \def\bacp  {Ba\-nach\wB Cen\-ter\wB Publ.}
   \def\bima  {Bulletin of the Institute of Mathematics}
   \def\coma  {Con\-temp.\wb Math.}
   \def\comp  {Com\-mun.\wb Math.\wb Phys.}
   \def\cpma  {Com\-pos.\wb Math.}
   \def\imrn  {Int.\wb Math.\wb Res.\wb Notices}
   \def\jlms  {J.\wB London\wB Math.\wb Soc.}
   \def\joal  {J.\wB Al\-ge\-bra}
   \def\jomp  {J.\wb Math.\wb Phys.}
   \def\jpaa  {J.\wB Pure\wB Appl.\wb Alg.}
   \def\mama  {Manuscripta\wB math.}
   \def\mams  {Memoirs\wB Amer.\wb Math.\wb Soc.}
   \def\nupb  {Nucl.\wb Phys.\ B}
   \def\slnm  {Sprin\-ger\wB Lecture\wB Notes\wB in\wB Mathematics}
   \def\taac  {Theo\-ry\wB and\wB Appl.\wb Cat.}

%%%%%%%%%%%%%%%%%%%%%%%%%%%%%%%%%%%%%%%%%%%%%%%%%%%%%%%%%%%%%%%%%%%%%%%%
\end{document}